\documentclass{article}

\usepackage[latin1]{inputenc}
\usepackage[T1]{fontenc}
\usepackage{amsmath}
\usepackage{amssymb}
\usepackage{theorem}
\usepackage{amscd}

\newtheorem{thm}{Theorem}
\newtheorem{defn}{Definition}
\newtheorem{lemma}{Lemma}

\newtheorem{exam}{Example}
\newtheorem{cory}{Corollary}
\newtheorem{rmk}{Remark}

\newcommand{\oct}{\mathbb{O}}
\newcommand{\R}{\mathbb{R}}
\newcommand{\C}{\mathbb{C}}
\newcommand{\h}{\mathbb{H}}
\newcommand{\g}{\mathfrak{g}}
\newcommand{\f}{\mathcal{F}}
\newcommand{\D}{\mbox{$D$\hspace{-.65em}/}}
\newcommand{\s}{\mathbb{R}^{2|2}}
\newcommand{\gt}{\tilde{\mathfrak{g}}}

\newcommand{\dl}[2]{\frac{\partial{#1}}{\partial{#2}}}
\newcommand{\te}[1]{\theta_{#1}}
\newcommand{\ta}{\theta}
\newcommand{\tb}{\bar{\theta}}
\newcommand{\lm}{\lambda}
\newcommand{\Lm}{\Lambda}
\newcommand{\iif}{if and only if }
\newcommand{\mnr}[1]{\mathfrak{M}_{#1}(\R)}
\newcommand{\mnc}[1]{\mathfrak{M}_{#1}(\C)}
\newcommand{\Gt}{G_{\tau}}
\newcommand{\Gtc}{G_{\tau}^{\C}}
\newcommand{\gtc}{\g_{\tau}^{\C}}
\newcommand{\gtau}{\g_{\tau}}
\newcommand{\mk}{\mathfrak{m}}
\newcommand{\gs}{\g_{\sigma}}
\newcommand{\gsc}{\g_{\sigma}^{\C}}
\newcommand{\Gs}{G_{\sigma}}
\newcommand{\Gsc}{G_{\sigma}^{\C}}
\author{Idrisse Khemar}
\title{Supersymmetric Harmonic Maps into Symmetric Spaces}
\date{}
\begin{document}
\maketitle

\section*{Introduction}

In this paper we study supersymmetric harmonic maps from the point
of view of integrable system. It is well known that harmonic maps
from $\R^2$ into a symmetric space are solutions of a integrable
system (see \cite{{DPW},{8},{12},{H1},{H2}}). We show here that the
superharmonic maps from $\s$ into a symmetric space are solutions of
a integrable system, more precisely of a first elliptic integrable
system in the sense of C.L. Terng (see \cite{tern}) and that we have a Weierstrass-type
representation in terms of holomorphic potentials (as well as of
meromorphic potentials). In the end of the paper we show that
superprimitive maps from $\s$ into a 4-symmetric space give us, by
restriction to $\R^2$, solutions of the second elliptic system
associated to the previous 4-symmetric space. This leads us to
conjecture that any second elliptic system associated to a
4-symmetric space has a geometrical interpretation in terms of
surfaces with values in a symmetric spaces, (such that a certain  associated
map is harmonic) as this is the case for Hamiltonian stationary
Lagrangian surfaces in Hermitian symmetric spaces (see \cite{HR3})
 or for $\rho$\,-harmonic surfaces of $\oct$ (see \cite{ki}).\\
Our paper is organized as follows. In the first section,  we define
superfields $\Phi\colon\s\to M$ from $\s$ to a Riemannian manifold,
and component fields. Then we recall the functor of points approach to
supermanifolds, we define the writing of a superfield and study its behaviour
 when we embedd the manifold $M$ in a Euclidiean space
$\R^N$. Lastly, we recall the derivation on $\s$. In section~\ref{2.2} we
 introduce the supersymmetric Lagrangian on $\s$, define the
 supersymmetric maps and derive the Euler-Lagrange equations in
 terms of the component fields. Next, we study the case $M=S^n$
 : we write the Euler-Lagrange equations in this case and we derive from
 them the  superharmonic maps equation in this case. Then we
 introduce the superspace formulation of the Lagrangian and derive
 the superharmonic maps equation for the general case of a Riemannian
 manifold $M$. In section~\ref{2.3}, we introduce the lift of a superfield
 with values in a symmetric space, then we express the superharmonic
 maps equation in terms of the Maurer-Cartan form of the lift. Once
 more, in order to make the comprehension easier, we first treat
 the case $M=S^n$, before the general case. In section~\ref{2.4}, we study
 the zero curvature equation (i.e. the Maurer-Cartan equation)
 for a 1-form on $\s$ with values in a Lie algebra. This allows to
 formulate the superharmonic maps equation  as the zero curvature
 equation for a 1-form on $\s$ with values in a loop space $\Lm\gtau$.
Then we precise the extended Maurer-Cartan form, and characterize
the superharmonic maps in terms of extended lifts. The section~\ref{2.5}
deals with the Weierstrass representation: we define holomorphic functions
and 1-forms in $\s$, and  then we define holomorphic potentials. We
show that we have a Weierstrass-type representation of the
superharmonic maps in terms of holomorphic potentials. Lastly, we
deal with meromorphic potentials. In section~\ref{weierfield}, we precise the
Weierstrass representation in terms of the component fields. In
section~\ref{2.7}, we study the superprimitive maps with values in a
4-symmetric spaces, and we precise their Weierstrass
representation. This allows us in the last section to show that
the restrictions to $\R^2$ of  superprimitive maps are
solutions of a second elliptic integrable system in the even part of
a super Lie algebra.

\section{Definitions and Notations}\label{2.1}

We consider the superspace $\s$ with coordinates $(x,y,\te{1},\te{2})$;
$(x,y)$ are the even coordinates and $(\te{1},\te{2})$ the odd
coordinates. Let $M$ be a Riemannian manifold. We will be interested
in maps $\Phi\colon\s\to M$ (which are even) i.e. morphisms of
sheaves of super $\R -$algebras from $\s$ to $M$ (see
\cite{{2},{berezin},{leites},{manin}}).
 We call these maps \emph{superfields}. We write such a superfield:
\begin{equation}\label{phi}
 \Phi=u + \te{1}\psi_1 + \te{2}\psi_2 + \te{1}\te{2}F'
\end{equation}
$u,\psi_1,\psi_2,F'$ are the component fields (see \cite{3}). We
view these as maps from $\R^2$ into a supermanifold: $u$ is a map
from $\R^2$ to $M$, $\psi_1,\psi_2$ are odd sections of $u^*(TM)$ and
$F'$ is a even section of $u^*(TM)$. So $u,F'$ are even whereas
$\psi_1,\psi_2$ are odd. The supermanifold of superfields $\Phi$ is
isomorphic  to the supermanifold of component fields
$\{u,\psi_1,\psi_2,F'\}$ (see \cite{3}). Besides the component
fields can be defined as the restriction to $\R^2$ of certain
derivatives of $\Phi$:
\begin{eqnarray}
u & = & i^*\Phi\colon\R^2\to M \nonumber\\
\psi_a & = & i^*D_a\Phi\colon \R^2\to u^*(\Pi TM)\label{i*}\\
F' & = & i^*(-\frac{1}{2}\varepsilon^{ab} D_aD_b \Phi)\colon \R^2\to
u^*(TM)\nonumber\ .
\end{eqnarray}
where $i\colon\R^2\to\s$ is the natural inclusion, $\Pi$ is the
functor which reverses the parity, and the left-invariant  vector
fields $D_a$ are defined below. This is the definition of the
component fields used in \cite{3}. We use another definition  based
on the morphism interpretation of superfields, which is equivalent
to the previous one, given by (\ref{i*}). Moreover as in
\cite{3} we use the functor of points approach to supermanifolds
(see \cite{2}). If $B$ is a supermanifold, then a $B-$point of $\s$
is a morphism $B\to\s$. It can be viewed as a family of points of
$\s$ parametrized by $B$, i.e. a section of the projection $\s\times
B\to B$. Then a map $\Phi$ from $\s$ to $M$ is a functor from the categoy of
supermanifolds, which to
each $B$ associates a map $\Phi_B\colon\s(B)\to M(B)$ from the set of
$B-$points of $\s$ to the set $M(B)$ of $B-$points of $M$. For example,
 if we take $B=\R^{0|L}$, which is the topogical space $\R^0$
endowed with the Grassman algebra $B_L=\R[\eta_1,\ldots,\eta_L]$ over
$\R^{L}$, then a $\R^{0|L}-$point of $\s$ is in the form $(x,y,\te{1},\te{2})$
where $x,y\in B_L^0$, the even part of $B_L$, and $\te1 ,\te2\in
B_L^1$, the odd part of $B_L$. Hence the set of $\R^{0|L}-$points
of $\s$ is $B_L^{2|2}:=(B_L^0)^2\times(B_L^1)^2$. Thus if we
restrict ourself to the category of supermanifolds $\R^{0|L},\,L\in
\mathbb{N}$, then a map $\Phi\colon\s\to M$ is a sequence $(\Phi_L)$, of
 $G^{\infty}$ functions defined by Rogers (\cite{rogers}), such that $\Phi_L$
is a $G^{\infty}$ function  from $B_L^{2|2}$ to the $G^{\infty}$ supermanifold
over $B_L$, $M(\R^{0|L})$, and such that ${\Phi_{L'}}_{|B_L^{2|2}}=\Phi_L$,
if $L \leq L'$.
Hence, in this case, if we suppose $M=\R^n$, we have
$M(\R^{0|L})=B_L^{n|0}=(B_L^0)^n$ and the writing (\ref{phi})
is the $z$ expansion of $\Phi_L$ (see \cite{rogers}). Further following
\cite{fp2}, we can say equivalently that if we denote by $\f$ the
infinite dimensional supermanifold of morphisms: $\s\to M$, then the
 functor defined by $\Phi$ is a functor
$B\mapsto\text{Hom}(B,\f)$: to each $B$ corresponds a $B-$point of
$\f$, i.e. a morphism $\Phi_B\colon\s\times B\to M$. It means that the map
$\Phi$ is a functor  which to each $B$ associates
a morphism of algebras $\Phi_B^*\colon C^{\infty}(M)\to C^{\infty}
(\s\times B)$.
In concrete terms, in all the paper, when we say: ``Let $\Phi\colon\s\to M$  be a map'',
one can consider that it means ``Let $B$ be a supermanifold and let
$\Phi_B\colon\s\times B\to M$ be a morphism'' (omitting the additional
 condition that $B\mapsto\Phi_B$ is functorial in $B$). $B$ can be viewed
as a ``space of parameters'', and $\Phi_B$ as a family of maps: $\s\to M$,
parametrized by $B$. We will never mention $B$ though it is tacitly assumed
to always be there. Moreover, when we speak about morphisms, these are even morphisms,
i.e. which preserve the parity, that is to say morphisms  of
 super $\R -$algebras. Thus as said above, a superfield is even.
But we will also be led to consider odd maps $A\colon\s\to M$, these are maps
 which give morphisms that reverse the parity.\\
 Let us now precise  the writing (\ref{phi}) and give our definition
 of the component fields.\\
 In the general case ($M$ is not an Euclidiean
 space $\R^N$) the formal writing (\ref{phi}) does not permit to have
 directly the morphism of super $\R-$algebras $\Phi^*$ as it happens in the
 case  $M=\R^N$, where the meaning of the writing (\ref{phi}) is clear:
 it is the writing of the morphism $\Phi^*$.
 Indeed, if $M=\R^N$ we have
 \begin{eqnarray}
 \forall f\in  C^{\infty}(\R^N), & &  \nonumber\\
  \Phi^*(f)=f\circ\Phi & = & f(u) + \sum_{k=1}^{\infty}\frac{f^{(k)}(u)}{k!}
\cdot(\te1\psi_1 +\te2\psi_2 + \te1\te2F')^k \nonumber\\
   & = & f(u) + \sum_{k=1}^{2}\frac{f^{(k)}(u)}{k!}
\cdot(\te1\psi_1 +\te2\psi_2 + \te1\te2 F')^k \nonumber\\
   & = & f(u) + \te1 df(u).\psi_1 +\te2 df(u).\psi_2  \nonumber\\
  &   & \qquad +\, \te1\te2 (df(u).F' -d^2f(u)(\psi_1,\psi_2))\label{phi*}
 \end{eqnarray}
 (we have used the fact that $\psi_1,\,\psi_2$ are odd).
Then we define the component fields as the the coefficient maps $a_I$
in the decomposition $\Phi=\sum\ta^I a_I$ in the morphism writing,
and as we will see below the equations (\ref{i*}) follow from this
definition.\\
In the general case, we must use local coordinates in $M$, to write
the morphism of algebras $\Phi^*$ in the same way as (\ref{phi*})
(see \cite{{berezin},{leites},{manin}}). But the
coefficient maps which appear in each chart in the equations
(\ref{phi*}) written in each chart, do not transform, through a change
of chart, in such a way that they define some unique
functions $u,\,\psi,\,F'$, which would allow us to give a sense to
(\ref{phi}) (in fact the coefficients corresponding to $u,\,\psi$
tranform correctely but not the one corresponding to $F'$). So the
writing (\ref{phi}) does not have any sense if we do not precise it. We
will  do it now. To do this we use the metric of $M$, more precisely
its Levi-Civita connection (it was already used  in the equation
(\ref{i*}), taken in \cite{3} as definition of  the component
fields, where the outer (leftmost) derivative in the expression of
$F'$ is a covariant derivative). We will show that for any
$\Phi\colon\s\to M$ there exist $u,\psi,F'$ which satisfy the
hypothesis above ($u,F'$ even, $\psi$ odd and $\psi,F'$ are tangent)
such that
\begin{eqnarray}
\forall f \in C^{\infty}(M), & & \nonumber\\
\Phi^*(f) & = & f(u) + \te1 df(u).\psi_1 +\te2 df(u).\psi_2\nonumber\\
 & & \qquad\quad +\ \te1\te2 (df(u).F' -(\nabla df)(u)(\psi_1,\psi_2))\label{phi*m}
\end{eqnarray}
where $\nabla df$ is the covariant derivative of $df$ (i.e. the
covariant Hessian of $f$): $(\nabla df)(X,Y)=\langle \nabla_X(\nabla
f),Y\rangle=\langle X,\nabla_Y(\nabla f)\rangle$. First, we remark
that if (\ref{phi*m}) is true, then $u,\psi,F'$ are unique. Then we
can define the component fields as being $u,\psi,F'$; and
(\ref{phi}) have a sense: it means that the morphism $\Phi^*$ is
given by (\ref{phi*m}).\\[.15cm]
 \indent Now, to prove (\ref{phi*m}), let us embedd isometrically $M$ in an
 Euclidiean space $\R^N$. Suppose first that $M$ is defined
 by a implicit equation in $\R^N$: $f(x)=0$, with
 $f\colon\R^N\to\R^{N-n}$ ($n=\dim M$). Then we have an isomorphism
 between $\{\text{superfields }\Phi\colon\s\to M\}$ and
 $\{\text{superfields }\Phi'\colon\s\to\R^N/{\Phi'}^*(f)=0\}$, the
 isomorphism is
\begin{equation}\label{iso}
  \Phi\longmapsto\Phi'=j\circ\Phi=(\,g\in C^{\infty}(\R^N)\mapsto
  \Phi^*(g_{|M})\,)
\end{equation}
where $j\colon M\to \R^N$ is the natural inclusion. In particular,
a superfield $\Phi'\colon\s\to\R^N$ is a superfield $\Phi$ from $\s$
into $M$ if and only if ${\Phi'}^*(f)=f\circ\Phi'=0$. It means that
if we write $\Phi'=u +\te1\psi_1 +\te2\psi_2 +\te1\te2 F$ then we
have by (\ref{phi*})
$$
0=f(u) + \te1 df(u).\psi_1 +\te2 df(u).\psi_2
 + \te1\te2 (df(u).F -d^2f(u)(\psi_1,\psi_2))
$$
hence $f(u)=0,\,df(u).\psi_a=0,\,df(u).F=d^2f(u)(\psi_1,\psi_2)$
\ i.e.
\begin{equation}\label{psitangent}
 \left\{ \begin{array}{l}
    u \text{ takes values in }M \\
    \psi_a \text{ takes values in }u^*(TM)\\
    df(u).F=d^2f(u)(\psi_1,\psi_2)\ .
     \end{array}\right.
\end{equation}
 Thus a superfield $\Phi'\colon\s\to\R^N$ is ``with values'' in $M$
if and only if $\Phi'=u +\te1\psi_1 +\te2\psi_2 +\te1\te2 F$ with
$(u,\psi,F)$ satisfying (\ref{psitangent}).\\
In the general case, there exists a family $(U_{\alpha})$ of open
sets in $\R^N$ such that $M\subset \bigcup_{\alpha} U_{\alpha}$ and
$C^{\infty}$ functions $f_{\alpha}\colon U_{\alpha}\to\R^{N-n}$ such
that $M\cap U_{\alpha}=f_{\alpha}^{-1}({0})$. Then
 $\Phi\mapsto j\circ\Phi$ is a isomorphism between $\{\Phi\colon\s\to M\}$
 and $\{\Phi'\colon\s\to\R^N/$ ${\Phi'}^*(f_{\alpha})=0,\,\forall
 \alpha\}$. When we write ${\Phi'}^*(f_{\alpha})=0$, it means that
 we consider $V_{\alpha}={\Phi'}^{-1}(U_{\alpha})$ (it is the open
 submanifold of $\s$ associated to $u^{-1}(U_{\alpha})\subset\R^2$,
 i.e. $u^{-1}(U_{\alpha})$ endowed with the restriction to
 $u^{-1}(U_{\alpha})$ of the structural sheaf of $\s$) and that
 $(\Phi_{|V_{\alpha}}')^*(f_{\alpha})=f_{\alpha}\circ\Phi_{|V_{\alpha}}'=0$.
 (see \cite{2}.) Hence a superfield $\Phi'\colon\s\to \R^N$ is with
 values in $M$ \iif  $\Phi'=u +\te1\psi_1 +\te2\psi_2 +\te1\te2 F$ with
$(u,\psi,F)$ satisfying (\ref{psitangent}) for each $f_{\alpha}$.
Now , we write that we have $\Phi^*(g_{|M})={\Phi'}^*(g)$, $\forall
g\in C^{\infty}(\R^N)$ :
$$
\Phi^*(g_{|M})=g(u) + \te1dg(u).\psi_1 + \te2dg(u).\psi_2 +
\te1\te2(dg(u).F-d^2g(u)(\psi_1,\psi_2)).
$$
Let $\mathrm{pr}(x)\colon\R^N\to T_xM$ be the orthogonal projection on $T_xM$
for $x\in M$, and $\mathrm{pr}^{\bot}(x)=Id-\mathrm{pr}(x)$; then set $F'=\mathrm{pr}(u).F$,
$F^{\bot}=\mathrm{pr}^{\bot}(u).F$, so that $F=F'+F^{\bot}$. Let also
$(e_1,\ldots,e_{N-n})$ be a local moving frame of $TM^{\bot}$. Then
we have
$$
dg(u).F-d^2g(u)(\psi_1,\psi_2)=\langle\nabla(g_{|M})(u),F'\rangle
+\langle \nabla g(u),F^{\bot}\rangle -\langle D_{\psi_1}\nabla g(u),\psi_2
 \rangle$$
 (where $D_{\psi_1}=\iota(\psi_1)d\,$). Now using that $\psi_1,\psi_2$
 are tangent to $M$ at $u$
 \begin{eqnarray*}
\langle D_{\psi_1}\nabla g(u),\psi_2 \rangle & = & \langle \mathrm{pr}(u).(D_{\psi_1}
\nabla g(u)),\psi_2 \rangle \\
 & = & \langle \mathrm{pr}(u).\left[D_{\psi_1}\left(\mathrm{pr}(\,).\nabla g\right)(u)+
  D_{\psi_1}\left(\mathrm{pr}^{\bot}(\,).\nabla g\right)(u)\right],\psi_2 \rangle\\
 & = & \langle \nabla_{\psi_1}\nabla (g_{|M}),\psi_2 \rangle +
\left\langle \mathrm{pr}(u).\left(D_{\psi_1}\sum_{i=1}^{N-n}\langle \nabla g,e_i \rangle
e_i\right),\psi_2\right\rangle\\
 & = & \nabla d(g_{|M})(u)(\psi_1,\psi_2) + \sum_{i=1}^{N-n}\langle
 \nabla g(u),e_i \rangle\langle de_i(u).\psi_1,\psi_2 \rangle
 \end{eqnarray*}
then
 \begin{eqnarray*}
dg(u).F-d^2g(u)(\psi_1,\psi_2) & = &  d(g_{|M})(u).F'-\nabla d(g_{|M})(u)
(\psi_1,\psi_2)\\
 & & +\,\langle \mathrm{pr}^{\bot}(u).\nabla g(u),F^{\bot}-\sum_{i=1}^{N-n}
 \langle de_i(u).\psi_1,\psi_2 \rangle e_i\rangle.
 \end{eqnarray*}
But, as $\Phi^*(g_{|M})$ depends only on $h=g_{|M}\in C^{\infty}(M)$,
we have
\begin{equation}\label{f''}
F^{\bot}=\sum_{i=1}^{N-n} \langle de_i(u).\psi_1,\psi_2 \rangle e_i
\end{equation}
and finally we obtain
\begin{eqnarray}
\forall h \in C^{\infty}(M), & & \nonumber\\
\Phi^*(h) & = & h(u) + \te1 dh(u).\psi_1 +\te2 dh(u).\psi_2\nonumber\\
 & & \qquad\quad +\ \te1\te2 (dh(u).F' -(\nabla dh)(u)(\psi_1,\psi_2))
 \label{8}
\end{eqnarray}
which is (\ref{phi*m}). And we have remarked that the coefficient maps
$\{u,\psi,F'\}$ are unique, so in particular they do not depend on
the embedding $M\hookrightarrow \R^N$. So we can define the
multiplet of the component fields of $\Phi$ in the general case: it is the
multiplet $\{u,\psi,F'\}$ which is defined by (\ref{phi*m}). It is an
intrinsec definition.\\
The isomorphism (\ref{iso}) leads to a isomorphim between the
component fields
$$\{u,\psi,F'\}\longmapsto \{u,\psi,F\}.$$
The only change is in the third component field. We have $F'=\mathrm{pr}(u).F$,
 and the orthogonal component $F^{\perp}$ of $F$ can be expressed in terms of
 $(u,\psi)$ as we can see it on (\ref{f''}) or on (\ref{psitangent}).\\
In the following when we consider a manifold $M$ with a natural
embedding $M\hookrightarrow \R^N$, we will identify $\Phi$ and
$\Phi'$, and we will talk about the two writings of $\Phi$: its
writing in $M$ and its writing in $\R^N$. But when we refer to
the component fields it will be always in $M$: $\{u,\psi,F'\}$. We
will in fact use only the writing in $\R^N$ because it is more
convenient to do computations, for example computations of
derivatives or multiplication of two superfields with values in a
Lie group, and because the meaning of the writing (\ref{phi}) in
$\R^N$ is clear and well known as well as how to  use it to do
computations. So we will not use the writing in $M$. Our aim was,
 first, to show that it is possible to generalize the writing
(\ref{phi}) in the general case of a Riemannian  manifold, then to
give a definition of the component fields which did not use the
derivatives of $\Phi$ (as in (\ref{i*})), and above all to show how
to deduce the component fields of $\Phi$ from its writing in $\R^N$: $u,\psi$
are the same and $F'=\mathrm{pr}(u).F$.
\begin{exam} $M=S^n\subset \R^{n+1}$.\\
A superfield $\Phi\colon\s\to\R^{n+1}$ is a superfield
$\Phi\colon\s\to S^n$ \iif \hbox{$\Phi^*(|\cdot|^2 -1)=(|\cdot|^2
-1)\circ\Phi=0$} ($|\cdot|$ being the Euclidiean norm in $\R^{n+1}$).
It means that
$$
0=\langle\Phi,\Phi \rangle -1= |u|^2-1 +2\te1 \langle\psi_1,u\rangle
+2\te2 \langle\psi_2,u\rangle +2\te1\te2 (\langle F,u\rangle -\langle\psi_1
,\psi_2\rangle)
$$
Thus $\Phi\colon\s\to\R^{n+1}$ takes values in $S^n$ \iif
$$
 \left\{ \begin{array}{l}
   u \text{ takes values in } S^n\\
    \psi_a\text{ is tangent to } S^n \text{ at } u  \\
  \langle F,u\rangle=\langle \psi_1,\psi_2\rangle
  \end{array}\right.
$$
In particular, in the case  of $S^n$ we have
$$ F^{\perp}=\langle\psi_1,\psi_2\rangle u.$$
\end{exam}

\subsection*{Derivation on $\s$.}
Let us introduce the left-invariant vector fields of $\s$:
\begin{eqnarray*}
D_1 & = & \dl{}{\te1} -\te1 \dl{}{x} -\te2\dl{}{y}\\
D_2 & = & \dl{}{\te2} -\te1 \dl{}{y} +\te2\dl{}{x}
\end{eqnarray*}
These vectors fields induce odd  derivations acting on
superfields $D_a\Phi=\iota(D_a)d\Phi$. Consider the case of
superfields with values in $\R^N$. Write $\Phi=u +\te1\psi_1 +\te2\psi_2 +
\te1\te2 F$ a superfield $\Phi\colon\s\to\R^N$. Then we have
\begin{eqnarray}
D_1\Phi & = & \psi_1 -\te1 \dl{u}{x} +\te2 \left(F-\dl{u}{y}\right)
+\te1\te2(\D\psi)_1 \label{d1}\\
D_2\Phi & = & \psi_2 -\te1 \left(\dl{u}{y}+F\right) +\te2 \dl{u}{x}
+\te1\te2(\D\psi)_2 \label{d2}
\end{eqnarray}
where
$$
\D\psi=\begin{pmatrix}
\displaystyle \dl{\psi_1}{y} -\dl{\psi_2}{x} \\
 \displaystyle -\dl{\psi_1}{x}-\dl{\psi_2}{y}
\end{pmatrix}=  \begin{pmatrix}
 \displaystyle   \dl{}{y} & \displaystyle -\dl{}{x}\\
 \displaystyle   -\dl{}{x} & \displaystyle -\dl{}{y}
  \end{pmatrix}
  \begin{pmatrix}
    \psi_1 \\
    \psi_2
  \end{pmatrix}
$$
Hence
$$
  \begin{array}{rclcrcl}
D_1D_1\Phi & = & -\displaystyle\dl{\Phi}{x} & , & D_1D_2\Phi & = &
 -R(\Phi)-\displaystyle\dl{\Phi}{y},\\
D_2D_1\Phi & = & R(\Phi)-\displaystyle\dl{\Phi}{y} & , & D_2D_2\Phi & = &
\displaystyle \dl{\Phi}{x},
  \end{array}
$$
where
\begin{eqnarray}
R(\Phi) & := & F+\te1 \left(\displaystyle \dl{\psi_2}{x}-\dl{\psi_1}{y}\right)
+ \te2 \left(\displaystyle \dl{\psi_1}{x}+\dl{\psi_2}{y}\right)
+\te1\te2 (\triangle u)\nonumber\\
       & := & F -\te1 (\D\psi)_1 - \te2 (\D\psi)_2
+\te1\te2 (\triangle u).\label{rphi}
\end{eqnarray}
Thus
\begin{eqnarray*}
D_1D_2-D_2D_1=-2R & , & [D_1,D_2]=D_1D_2 +D_2D_1=-2\dl{}{y} \\
{[D_1,D_1]}=2 D_{1}^{2}= -2\dl{}{x}  & , &  [D_2,D_2]= 2\dl{}{x}
\end{eqnarray*}
(In all the paper, we denote by $[\,,\,]$ the superbracket in the considered
super Lie algebra).\\
Let us set
\begin{eqnarray*}
D & = & \frac{1}{2}(D_1-iD_2)=\dl{}{\ta} -\ta\dl{}{z} \\
\bar{D} & = & \frac{1}{2}(D_1+iD_2)=\dl{}{\tb} -\tb\dl{}{\bar{z}}
\end{eqnarray*}
where $\ta=\te1 +i\te2$, $\dl{}{\ta}=\frac{1}{2}(\dl{}{\te1}
-i\dl{}{\te2})$. Setting $\psi=\psi_1 -i\psi_2$, we can write $\Phi=u
+ \frac{1}{2}(\ta\psi +\tb\bar{\psi})+\frac{i}{2}\ta\tb F$,
thus
\begin{eqnarray}
D\Phi & = & \frac{1}{2}\psi -\ta\dl{u}{z} +\frac{i}{2}\tb F
-\frac{1}{2}\ta\tb\dl{\bar{\psi}}{z}\label{d}\\
\bar{D}\Phi & = & \frac{1}{2}\bar{\psi} -\tb\dl{u}{\bar{z}} -\frac{i}{2}\ta F
+\frac{1}{2}\ta\tb\dl{\psi}{\bar{z}}\label{dbar}
\end{eqnarray}
Then
$$ \begin{array}{rclll}
D\bar{D} & = & \displaystyle\frac{1}{4}(D_1 -iD_2)(D_1 + iD_2) & = &
\displaystyle\frac{1}{4}(D_1^2 +D_2^2 + i(D_1D_2 -D_2D_1))\\
  & = & \displaystyle\frac{i}{4}(D_1D_2 -D_2D_1) & = &
  \displaystyle -\frac{i}{2}R
  \end{array}
$$
hence
$$ D\bar{D}=-\bar{D}D=-\frac{i}{2}R.$$
We have also $D^2=-\dl{}{z}$, $\bar{D}^2=-\dl{}{\bar{z}}$. Let us
compute $\bar{D}D\Phi$:
\begin{eqnarray}
\bar{D}D\Phi & = & \bar{D}\left(\frac{1}{2}\psi -\ta\dl{u}{z} +\frac{i}{2}\tb F
-\frac{1}{2}\ta\tb\dl{\bar{\psi}}{z}\right)\nonumber\\
 & = & \frac{i}{2}F +\frac{\ta}{2}\dl{\bar{\psi}}{z}
 -\frac{\tb}{2}\dl{\psi}{\bar{z}}
 -\ta\tb\dl{}{\bar{z}}\left(\dl{u}{z}\right)\nonumber\\
  & = & \frac{i}{2}F + i\,\text{Im}\left(\ta\dl{\bar{\psi}}{z}\right)
  -\frac{\ta\tb}{4}(\triangle u)\label{rddphi}.
\end{eqnarray}
Let us denote by $i\colon\R^2\to\s$ the natural inclusion, then using
(\ref{d1})-(\ref{d2}) and (\ref{rphi}) we have
\begin{eqnarray*}
u & = & i^*\Phi\\
\psi_a & = & i^*D_a\Phi\\
F & = & i^*(-\frac{1}{2}\varepsilon^{ab} D_aD_b \Phi)
\end{eqnarray*}
and we recover (\ref{i*}) for $M=\R^N$.\\
Let us return to the general case of superfields with values in
$M$. In order to write (\ref{i*}) in $M$, we need a covariant
derivative in the expression of $F'$ to define the action of $D_a$
on a section of the bundle $\Phi^* TM$. In order to do this we use
the pullback of the Levi-Civita connection. Suppose that $M$ is
isometrically embedded in $\R^N$. Let $X$ be a section of $\Phi^*
TM$ (for example $X=D_b\Phi$) then using the writing in $\R^N$ (i.e.
considering that a map with values in $M$ takes values in $\R^N$)
we have
$$
\nabla_{D_a}X=\mathrm{pr}(\Phi).D_a X\ .
$$
Let us precise the expression $\mathrm{pr}(\Phi).D_a X$. The projection
$\mathrm{pr}$ is a map from $M$ into $\mathcal{L}(\R^N)$, the
algebra of endomorphisms of $\R^N$. We consider $\mathrm{pr}\circ\Phi$ which
we write $\mathrm{pr}(\Phi)$. Then considering the maps $\mathrm{pr}(\Phi)\colon\s\to
\mathcal{L}(\R^N)$, $D_a X\colon\s\to\R^N$, and $B\colon (A,v)\in \mathcal{L}
(\R^N)\times\R^N\mapsto A.v$, we form $B(\mathrm{pr}(\Phi),D_a
X)\colon\s\to\R^N$. Now, since $\mathcal{L}(\R^N)$ is a finite
dimensional vector space we can write from (\ref{phi*m}):
\begin{eqnarray*}
\mathrm{pr}(\Phi) = \Phi^*(\mathrm{pr}) & = & \mathrm{pr}(u) + \te1
d\mathrm{pr}(u).\psi_1 +\te2 d\mathrm{pr}(u).\psi_2\\
 & & \qquad\quad +\ \te1\te2 (d\mathrm{pr}(u).F' -(\nabla d\mathrm{pr})(u)(\psi_1,\psi_2))
\end{eqnarray*}
(we can not use (\ref{phi*}) because $\mathrm{pr}$ is only defined on $M$). This
is the writing of the superfield $\mathrm{pr}\circ\Phi\colon\s\to\mathcal{L}(\R^N)$,
so we can write
$$i^*(\nabla_{D_a}D_b\Phi)=i^*(\mathrm{pr}(\Phi).D_a D_b\Phi)=\mathrm{pr}(u).i^*(D_a
D_b\Phi)$$
 thus $i^*(-\frac{1}{2}\varepsilon^{ab}\nabla_{D_a}
D_b\Phi)=\mathrm{pr}(u).F=F'.$ So we have (\ref{i*}) in the general case.

\begin{exam} $M=S^n\subset \R^{n+1}$.\\
We have $\mathrm{pr}(x)=Id - \langle\cdot,x\rangle x$ for $x\in S^n$. So for
$X$ a section of $\Phi^* TS^n$, we have
$$\nabla_{D_a}X=D_a X - \langle D_a X,\Phi\rangle\Phi.$$
\end{exam}

\section{Supersymmetric Lagrangian}\label{2.2}

\subsection{Euler-Lagrange equations}

We consider the following supersymmetric Lagrangian (see \cite{3}):
\begin{equation}\label{lagr}
L=-\frac{1}{2}|du|^2 +\frac{1}{2}\langle\psi\D_u\psi\rangle
+\frac{1}{12}\varepsilon^{ab}\varepsilon^{cd}\langle\psi_a,R(\psi_b,\psi_c)
\psi_d\rangle +\frac{1}{2}|F'|^2
\end{equation}
where $\langle\psi\D_u\psi\rangle=\langle\psi_1,(\D_u\psi)_2\rangle-
\langle\psi_2,(\D_u\psi)_1\rangle$, $R$ is the curvature of $M$ and
$$
\D_u\psi=\begin{pmatrix}
\displaystyle \dl{\psi_1}{y} -\dl{\psi_2}{x} \\
 \displaystyle -\dl{\psi_1}{x}-\dl{\psi_2}{y}
\end{pmatrix}
$$
(\,$\dl{\psi_k}{x_i}$ is of course a covariant derivative).
This Lagrangian can be obtained by reduction to $\s$ of the supersymmetric
$\sigma-$model Lagrangian on $\R^{3|2}$ (see \cite{3}).
We associate to this Lagrangian the action
$\mathcal{A}(\Phi)=\int L(\Phi)dxdy$. It is a functional on the
multiplets of components fields $\{u,\psi,F'\}$ of superfields
 $\Phi\colon\s\to M$, which is supersymmetric.
 \begin{defn}
A superfield $\Phi\colon\s\to M$ is superharmonic if it is a
critical point of the action $\mathcal{A}$
\end{defn}

\begin{thm}
 If we suppose that $\nabla R=0$ in $M$ (the covariant derivative of
 the curvature vanishes) then the Euler-Lagrange equations
 associated to the action $\mathcal{A}$ are:
 \begin{equation}\label{elm}
 \boxed{
 \begin{array}{rcl}
 \triangle u & = & \displaystyle\frac{1}{2}(R(\psi_1,\psi_1)-R(\psi_2,\psi_2))\dl{u}{x}
 + R(\psi_1,\psi_2)\dl{u}{y}\\
 \D_u\psi & = & \displaystyle
  \begin{pmatrix}
  R(\psi_1,\psi_2)\psi_1   \\
  -R(\psi_1,\psi_2)\psi_2
  \end{pmatrix}\\
  F' & = & 0
\end{array}}
\end{equation}
\end{thm}
\textbf{Proof.} We compute the variation  of each term in
the Lagrangian, keeping in mind that $\psi_1,\psi_2$ are odd (so their
coordinates anticommutate $\psi_1^{i}\psi_2^{j}=-\psi_2^{j}\psi_1^{i})$:\\

\noindent$\bullet\ \displaystyle\delta(\frac{1}{2}|du|^2)=\langle-
  \triangle u ,\delta u \rangle  + \text{div}(\langle du,\delta u\rangle)$

 \begin{eqnarray*}
\bullet\ \displaystyle\delta(\frac{1}{2}\langle\psi\D_u\psi\rangle) & = &
 \frac{1}{2}(  \langle\delta_{\nabla}\psi_1,(\D_u\psi)_2\rangle +
\langle\psi_1,\delta_{\nabla}(\D_u\psi)_2\rangle\hfill\null\\
 & &  - \langle\delta_{\nabla}\psi_2,(\D_u\psi)_1\rangle
-\langle\psi_2,\delta_{\nabla}(\D_u\psi)_1\rangle)\\
 & = & \frac{1}{2}\left[\frac{}{}\langle\delta_{\nabla}\psi_1,(\D_u\psi)_2\rangle
-\langle\delta_{\nabla}\psi_2,(\D_u\psi)_1\rangle\right.\\
 & &  + \displaystyle\left\langle\psi_1,-\dl{}{x}\delta_{\nabla}\psi_1 -\dl{}{y}
 \delta_{\nabla} \psi_2\right\rangle -\left\langle \psi_2,\dl{}{y}\delta_{\nabla}\psi_1
 -\dl{}{x} \delta_{\nabla}\psi_2\right\rangle\\
 & &  + \displaystyle\left\langle \psi_1,R\left(\delta u,-\dl{u}{x}\right)\psi_1
 -R\left(\delta u,-\dl{u}{y}\right)\psi_2\right\rangle\\
 & & -\left. \displaystyle\left\langle \psi_2,R\left(\delta u,\dl{u}{y}\right)\psi_1
 +R\left(\delta u,-\dl{u}{x}\right)\psi_2\right\rangle \right]
\end{eqnarray*}
we have used \,$\delta_{\nabla}\dl{\psi_k}{x_i}-\dl{}{x_i}\delta_{\nabla}\psi_k
=R(\delta u,\dl{u}{x_i})\psi_k$. Then we write that
\begin{eqnarray*}
\left\langle\psi_a,\dl{}{x_i}\delta_{\nabla}\psi_b\right\rangle & = &
  -\left\langle\dl{\psi_a}{x_i},\delta_{\nabla}\psi_b\right\rangle +
\dl{}{x_i}\left\langle \psi_a,\delta_{\nabla}\psi_b\right\rangle \\
 & = & \left\langle \delta_{\nabla}\psi_b,\dl{\psi_a}{x_i}\right\rangle
+ \dl{}{x_i}\left\langle \psi_a,\delta_{\nabla}\psi_b\right\rangle
\end{eqnarray*}
and that
$$
\left\langle\psi_a,R\left(\delta u,\dl{u}{x_i}\right)\psi_b\right\rangle =
\left\langle\ R\left(\psi_b,\psi_a \right)\dl{u}{x_i},\delta u
\right\rangle
$$
thus we obtain
\begin{multline*}
\displaystyle\delta(\frac{1}{2}\langle\psi\D_u\psi\rangle) =\\
\frac{1}{2}\left[\left\langle\delta_{\nabla}\psi_1,(\D_u\psi)_2 + \left(
-\dl{\psi_1}{x} -\dl{\psi_2}{y}\right)\right\rangle
 -\left\langle\delta_{\nabla}\psi_2,(\D_u\psi)_1 + \left(\dl{\psi_1}{y}
 -\dl{\psi_2}{x}\right)\right\rangle\right.\\
 + \dl{}{x}\left(-\langle\psi_1,\delta_{\nabla}\psi_1\rangle +
 \langle\psi_2,\delta_{\nabla}\psi_2\rangle\right) +
 \dl{}{y}\left(-\langle\psi_1,\delta_{\nabla}\psi_2\rangle -
 \langle\psi_2,\delta_{\nabla}\psi_1\rangle\right)\\
 -\left.\left\langle\left( R\left(\psi_1,\psi_1 \right)\dl{u}{x}
+ R\left(\psi_2,\psi_1 \right)\dl{u}{y} + R\left(\psi_1,\psi_2\right)\dl{u}{y}
- R\left(\psi_2,\psi_ 2\right)\dl{u}{x}\right),\delta u \right\rangle\right]
\end{multline*}
and finally
\begin{multline*}
\displaystyle\delta(\frac{1}{2}\langle\psi\D_u\psi\rangle) =
\langle\delta_{\nabla}\psi_1,(\D_u\psi)_2\rangle-
\langle\delta_{\nabla}\psi_2,(\D_u\psi)_1\rangle\\
-\left\langle \displaystyle\left[\frac{1}{2}(R(\psi_1,\psi_1)-R(\psi_2,\psi_2))
\dl{u}{x} + R(\psi_1,\psi_2)\dl{u}{y}\right],\delta u \right\rangle \\
+ \text{div}(\cdots)
\end{multline*}
$$
\begin{array}{l}\negthickspace\negthickspace\negthickspace
\bullet \ \displaystyle\delta\left(\frac{1}{12}\varepsilon^{ab}\varepsilon^{cd}
\langle\psi_a, R(\psi_b,\psi_c)\psi_d\rangle\right)\hfill\null\\
\begin{array}{rcl}
 & = & \displaystyle\frac{1}{12}\varepsilon^{ab}\varepsilon^{cd}\left(
\nabla_{\delta u} R(\psi_b,\psi_c,\psi_d,\psi_a) +
R(\delta\psi_a,\psi_b,\psi_c,\psi_d)\right. \\
 & & \hfill\left. +R(\psi_a,\delta\psi_b,\psi_c,\psi_d)
+R(\psi_a,\psi_b,\delta\psi_c,\psi_d) +R(\psi_a,\psi_b,\psi_c,\delta\psi_d)
\right)\hfill\null\\
 & = & \displaystyle\frac{1}{12}\varepsilon^{ab}\varepsilon^{cd}
(0+\langle\delta\psi_a,R(\psi_b,\psi_c)\psi_d\rangle
+\langle\delta\psi_b,R(\psi_d,\psi_a)\psi_c\rangle  \\
 & & \hfill +\displaystyle\frac{}{}\langle\delta\psi_c,R(\psi_d,\psi_a)\psi_b\rangle
 +\langle\delta\psi_d,R(\psi_b,\psi_c)\psi_a
 \rangle \hfill\null\\
 & & \hfill\displaystyle\frac{}{}(\text{using the symmetries of }R\,)\hfill\null\\
 & = & \displaystyle\frac{1}{12}(
\langle\delta\psi_1, R(\psi_2,\psi_1)\psi_2 - R(\psi_2,\psi_2)\psi_1
-R(\psi_2,\psi_2)\psi_1 + R(\psi_1,\psi_2)\psi_2 \\
 & & \hfill  +R(\psi_2,\psi_1)\psi_2 -R(\psi_2,\psi_2)\psi_1
  - R(\psi_2,\psi_2)\psi_1 + R(\psi_1,\psi_2)\psi_2 \rangle \hfill\null\\
 & &  + \langle\delta\psi_2,- R(\psi_1,\psi_1)\psi_2 + R(\psi_1,\psi_2)\psi_1
 + R(\psi_2,\psi_1)\psi_1 - R(\psi_1,\psi_1)\psi_2 \\
 & & \hfill  - R(\psi_1,\psi_1)\psi_2 + R(\psi_1,\psi_2)\psi_1
 + R(\psi_2,\psi_1)\psi_1 - R(\psi_1,\psi_1)\psi_2 \rangle ) \hfill\null\\
 & = & \displaystyle\frac{1}{12}(
\langle\delta\psi_1, -4R(\psi_2,\psi_2)\psi_1 + 4R(\psi_1,\psi_2)\psi_2
\rangle \\
 & & \hfill + \langle\delta\psi_2, 4R(\psi_2,\psi_1)\psi_1 -
4R(\psi_1,\psi_1)\psi_2\rangle)\hfill\null\\
 & = & \displaystyle\frac{1}{3}(
\langle\delta\psi_1, R(\psi_1,\psi_2)\psi_2 -
R(\psi_2,\psi_2)\psi_1\rangle\\
 & & \hfill + \langle\delta\psi_2, R(\psi_2,\psi_1)\psi_1 - R(\psi_1,\psi_1)\psi_2
 \rangle).\hfill\null
\end{array}
\end{array}
$$
Finally, by using the Bianchi identity we obtain:
$$\displaystyle\delta\left(\frac{1}{12}\varepsilon^{ab}\varepsilon^{cd}
\langle\psi_a, R(\psi_b,\psi_c)\psi_d\rangle\right)=
\langle \delta_{\nabla}\psi_1 ,R(\psi_1,\psi_2)\psi_2\rangle
+\langle\delta_{\nabla}\psi_2,R(\psi_2,\psi_1)\psi_1\rangle.$$

\noindent$\bullet \ \displaystyle\delta\left(\frac{1}{2}|F'|^2\right)=
\langle F',\delta_{\nabla}F'\rangle$\\

\noindent Hence the first variation of the Lagrangian is:
\begin{multline*}
\delta\mathcal{L}=\int\left[\left\langle \triangle u - \displaystyle\frac{1}{2}
(R(\psi_1,\psi_1)-R(\psi_2,\psi_2))\dl{u}{x} - R(\psi_1,\psi_2)\dl{u}{y}
,\delta u\right\rangle\right.\\
+\left\langle\delta_{\nabla}\psi_1,(\D_u\psi)_2 +R(\psi_1,\psi_2)\psi_2
\right\rangle
-\left\langle\delta_{\nabla}\psi_2,(\D_u\psi)_1 -R(\psi_1,\psi_2)\psi_1
\right\rangle\\
+\left. \langle F',\delta_{\nabla}F'\rangle\frac{}{}\right]dxdy
\end{multline*}
This completes the proof of the theorem.\hfill$\blacksquare$\\
\begin{rmk}\emph{
In any symmetric space, $\nabla R=0$, so that the preceding result
holds. Moreover in the general case of a Riemannian
manifold $M$ the Euler-Lagrange equations are obtained by adding to
the right hand side of the first equation of (\ref{elm}) the term
$-\frac{1}{2}(\nabla_{\psi_1}R)(\psi_1,\psi_2)\psi_2$.}
\end{rmk}
\subsection{The case $M=S^n$.}
The curvature of $S^n$ is given by
\begin{eqnarray*}
R(X,Y,Z,T) & = & \langle X,T\rangle\langle Y,Z\rangle -\langle X,Z \rangle
\langle T,Y\rangle\\
 & = & (\delta^{il}\delta^{jk} -\delta^{ik} \delta^{jl})X_i Y_j Z_k  T_l
\end{eqnarray*}
so
\begin{eqnarray*}
R(V_1,V_2)V_3 &=& \langle V_2,V_3\rangle V_1 + \langle V_1,V_3\rangle V_2\\
R(V_1,V_2)Z &=& -\langle V_2,Z \rangle V_1 -\langle V_1,Z\rangle V_2
\end{eqnarray*}
where $V_1,V_2,V_3$ are odd and $Z$ is even.\\
Thus the Euler-Lagrange equations for $S^n$ are :
\begin{eqnarray*}
\triangle u + |du|^2u & = & -\displaystyle\left\langle\psi_1,\dl{u}{x}
\right\rangle\psi_1
 +\left\langle\psi_2,\dl{u}{x}\right\rangle\psi_2\\
 & & \qquad - \left(\left\langle\psi_2,\dl{u}{y}\right\rangle\psi_1 +\left\langle\psi_1,
 \dl{u}{y}\right\rangle\psi_2\right) \\
 \D_u\psi & = &
  \begin{pmatrix}
  \langle\psi_2,\psi_1\rangle\psi_1   \\
  \langle\psi_2,\psi_1\rangle\psi_2
  \end{pmatrix}\\
  F & = & \langle\psi_1,\psi_2\rangle u
\end{eqnarray*}
Let us now rewrite these equations by using the complex variable and
setting $\psi =\psi_1 -i\psi_2$:
\begin{equation}\label{elsnc}
\boxed{
  \begin{array}{rcl}
4\displaystyle\dl{{}^{\nabla}}{\bar z }\left(\dl{u}{z}\right) & = &
\displaystyle\left(\psi\left\langle\psi ,\dl{u}{\bar{z}}\right\rangle
 + \bar\psi\left\langle\bar\psi ,\dl{u}{z}\right\rangle
 \right)\\[.3cm]
\displaystyle\dl{{}^{\nabla}\psi}{\bar{z}} & = & \displaystyle\frac{1}{4}
\langle\bar{\psi},\psi\rangle \bar{\psi}\\[.25cm]
 F &= & \displaystyle\frac{1}{2i}\langle \psi,\bar{\psi}\rangle u
\end{array}}
\end{equation}

\begin{thm}
Let $\Phi\colon\s\to S^n$ be a superfield, then $\Phi$ is
superharmonic \iif
\begin{equation}\label{dbard=0}
  \bar{D}D\Phi +\langle\bar{D}\Phi,D\Phi\rangle\Phi=0
\end{equation}
in $\R^{n+1}$.
\end{thm}
\textbf{Proof.} According to (\ref{rddphi}), we have
$$
\bar{D}D\Phi =\frac{i}{2}F + i\,\text{Im}\left(\ta\dl{\bar{\psi}}{z}\right)
-\ta\tb\dl{}{\bar{z}}\left(\dl{u}{z}\right).
$$
Moreover, by using (\ref{d}),(\ref{dbar})
\begin{multline*}
\langle\bar{D}\Phi,D\Phi\rangle\Phi=\frac{1}{4}\langle\bar{\psi},\psi\rangle
+\ta\left(\frac{1}{2}\left\langle\bar{\psi},\dl{u}{z} \right\rangle -
\frac{i}{4}\langle F,\psi \rangle\right)\\
 +\tb\left(-\frac{1}{2}\left\langle \dl{u}{\bar{z}},\psi\right\rangle
  -\frac{i}{4}\langle\bar{\psi},F \rangle\right)\\
+\ta\tb \left( -\frac{1}{4}\left\langle \bar{\psi},\dl{\bar{\psi}}{z}\right
\rangle +  \frac{1}{4}\left\langle \dl{\psi}{\bar{z}},\psi\right\rangle +
\frac{1}{4}|F|^2 -\left\langle \dl{u}{\bar{z}},\dl{u}{z}\right\rangle\right).
\end{multline*}
But since $\langle\psi, u\rangle=\langle\bar{\psi},u\rangle=0$ we have
$\langle\bar{\psi},\dl{u}{z}\rangle =-\langle
\dl{\bar{\psi}}{z},u\rangle$ and $\langle \dl{u}{\bar{z}},\psi\rangle =
-\langle u,\dl{\psi}{\bar{z}}\rangle$  so
\begin{multline*}
\langle\bar{D}\Phi,D\Phi\rangle\Phi=\frac{1}{4}\langle\bar{\psi},\psi\rangle
-\ta\left(\frac{1}{2}\left\langle\dl{\bar{\psi}}{z},u\right\rangle +
\frac{i}{4}\langle F,\psi \rangle\right)\\
 +\tb\left(\frac{1}{2}\left\langle \dl{\psi}{\bar{z}},u\right\rangle
  -\frac{i}{4}\langle\bar{\psi},F \rangle\right)\\
+\ta\tb \left( \frac{1}{2}\text{Re}\left(\left\langle \dl{\psi}{\bar{z}},\psi
\right\rangle\right) +
\frac{1}{4}|F|^2 -\left\langle \dl{u}{\bar{z}},\dl{u}{z}\right\rangle\right).
\end{multline*}
Hence
\begin{multline*}
\bar{D}D\Phi +\langle\bar{D}\Phi,D\Phi\rangle\Phi \\
 =\bar{D}D\Phi +\langle\bar{D}\Phi,D\Phi\rangle \left(u
+ \displaystyle\frac{1}{2}(\ta\psi +\tb\bar{\psi})+\frac{i}{2}\ta\tb F
\right)\qquad\qquad\qquad\\
  =  \displaystyle\left(\frac{i}{2}F + \frac{1}{4}\langle\bar{\psi},
 \psi\rangle u\right)\qquad\qquad\qquad\qquad\qquad\qquad\qquad\qquad\qquad
 \,\\
   +\frac{\ta}{2}\left(\dl{\bar{\psi}}{z} -\left\langle\dl{\bar{\psi}}{z}
 ,u\right\rangle u + \frac{1}{4}\langle\bar{\psi},\psi\rangle\psi
-\frac{i}{2}\langle F,\psi \rangle u\right)\\
   +\frac{\tb}{2}\left(-\dl{\psi}{\bar{z}} +\left\langle\dl{\psi}{\bar{z}}
 ,u\right\rangle u + \frac{1}{4}\langle\bar{\psi},\psi\rangle\bar{\psi}
-\frac{i}{2}\langle \bar{\psi},F \rangle u\right)\\
\qquad +\ta\tb\left(-\left[\dl{}{\bar{z}}\dl{u}{z} +\left\langle \dl{u}{\bar{z}}
 ,\dl{u}{z}\right\rangle\right] +\frac{1}{4}\left[\psi\left\langle\psi,
 \dl{u}{\bar{z}}\right\rangle +\bar{\psi}\left\langle \bar\psi ,
 \dl{u}{z}\right\rangle\right]\right. \\
 +\frac{i}{8}\langle  F,\psi\rangle\bar\psi -\frac{i}{8}\langle\bar\psi ,F
 \rangle\psi \\
 +\left[ \frac{1}{4}|F|^2 + \frac{1}{2}\text{Re}\left(\left\langle \dl{\psi}
 {\bar{z}},\psi\right\rangle\right)\right]u +\left.\frac{i}{8}\langle\bar\psi
 ,\psi\rangle F \right).
\end{multline*}
So we see that if $\Phi$ satisfies (\ref{elsnc}) then this
expression vanishes because $\langle F,\psi\rangle=\langle F,\bar\psi\rangle
=0$ and $\text{Re}\left(\left\langle \dl{\psi}
 {\bar{z}},\psi\right\rangle\right)=\text{Re}\langle\bar\psi
 ,\psi {\rangle}^2=-4|F|^2$ by using (\ref{elsnc}).\\
 Conversely, if this expression vanishes then the vanishing of the
 first term gives us the third equation of (\ref{elsnc}), thus we
 have $\langle F,\psi\rangle=0$ and so the vanishing of the therm in
 $\ta$ gives us the second equation of (\ref{elsnc}). Lastly the
 first equation of (\ref{elsnc}) is given by the vanishing of the
 term in $\ta\tb$ and by using the second and third equation of
 (\ref{elsnc}). This completes the proof.\hfill $\blacksquare$\\
\begin{rmk}\emph{
 The equation (\ref{dbard=0}) is the analogue of
  the  equation for harmonic maps $u\colon\R^2 \to S^n$:
  $$
  \dl{}{\bar{z}}\left(\dl{u}{z}\right) +\left\langle \dl{u}{\bar{z}} ,
  \dl{u}{z}\right \rangle=0.
  $$
In fact,  equation (\ref{dbard=0}) means that
$$
\nabla_{\bar{D}}D\Phi=0.
$$
Indeed we have
$\nabla_{\bar{D}}D\Phi=\mathrm{pr}(\Phi).\bar{D}D\Phi=\bar{D}D\Phi-\langle \bar{D}D\Phi
,\Phi\rangle\Phi$ but
$$
\begin{array}{rccll}
\langle \bar{D}D\Phi,\Phi\rangle\Phi &= & \bar{D}(\langle
D\Phi,\Phi\rangle)  & + & \langle D\Phi,\bar{D}\Phi\rangle \\
 & = &  0  & - & \langle \bar{D}\Phi,D\Phi\rangle
\end{array}
$$
because $\langle \Phi,\Phi\rangle=1 \Longrightarrow \langle D\Phi,\Phi\rangle
=0.$ So
$$
\nabla_{\bar{D}}D\Phi=\bar{D}D\Phi + \langle\bar{D}\Phi,D\Phi\rangle\Phi.
$$
It is a general result that $\Phi\colon\s\to M$ (Riemannian without
other hypothesis) is superharmonic \iif $\nabla_{\bar{D}}D\Phi=0$.
To prove it we need to use the superspace formulation  for the
supersymmetric Lagrangian. This is what we are going to do now.}
\end{rmk}

\subsection{The superspace formulation}

We consider the Lagrangian density  on $\s$ (see \cite{3}):
$$
L_0=dxdyd\te1 d\te2 \frac{1}{4}\varepsilon^{ab}\langle
D_a\Phi,D_b\Phi\rangle .
$$
$\Phi$ is a superfield $\Phi\colon\s\to M$, and $\langle\cdot,\cdot\rangle$
is the metric on $M$ pulled back to a metric on $\Phi^*TM$. Then,
according to \cite{3} the supersymmetric Lagrangian $L$, given in
(\ref{lagr}), is obtained by integrating over the $\ta$ variables
the Lagrangian density:
$$
L=\int d\te1 d\te2\frac{1}{4}\varepsilon^{ab}\langle D_a\Phi,D_b\Phi\rangle .
$$
Let us compute the variation of $L_0$ under an arbitrary even
variation $\delta\Phi$ of the superfield $\Phi$. We will set
$\nabla_{D_a}=D_a^{\nabla}$. Then, following \cite{3}, we have
\begin{eqnarray*}
\delta L_0 & = & dxdyd\te1 d\te2 \frac{1}{4}\varepsilon^{ab}(\langle
\delta_{\nabla}D_a\Phi,D_b\Phi\rangle + \langle D_a\Phi,\delta_{\nabla}
 D_b\Phi\rangle)\\
 & = & dxdyd\te1 d\te2 \frac{1}{2}\varepsilon^{ab}\langle \delta_{\nabla}
D_a\Phi,D_b\Phi\rangle\\
 & = & dxdyd\te1 d\te2 \frac{1}{2}\varepsilon^{ab}\langle D_a^{\nabla}
\delta_{\nabla}\Phi,D_b\Phi\rangle\\
& = & dxdyd\te1 d\te2 \frac{1}{2}\varepsilon^{ab}(D_a\langle \delta \Phi,
D_b\Phi\rangle -\langle \delta\Phi, D_a^{\nabla}D_b\Phi\rangle )\\
 & = & d\left[\iota(D_a)\left(dxdyd\te1 d\te2 \frac{1}{2}\varepsilon^{ab}\langle
D_b\Phi,\delta\Phi\rangle\right)\right] \\
 & &  \qquad\qquad\qquad - dxdyd\te1 d\te2 \frac{1}{2}\langle\delta\Phi,
(D_1^{\nabla}D_2 - D_2^{\nabla}D_1)\Phi\rangle.
\end{eqnarray*}
we have used  at the last stage the fact that the density $dxdyd\te1
d\te2$ is invariant under $D_a$ and the Cartan formula for the Lie
derivative. So the Euler-Lagrange equation in superspace is
$$
(D_1^{\nabla}D_2 - D_2^{\nabla}D_1)\Phi=0
$$
or equivalently,
\begin{equation}\label{dcovd=0}
\bar{D}^{\nabla}D\Phi=0
\end{equation}

\section{Lift of a superharmonic map into a symmetric space}\label{2.3}

\subsection{The case $M=S^n$}

We consider the quotient map $\pi\colon \text{SO}(n+1)\to S^n$ defined by
$\pi(v_1,\ldots,v_{n+1})=v_{n+1}$. We will say that
$\f\colon\s\to\text{SO}(n+1)$ is a lift of $\Phi\colon\s\to S^n$ if
$\pi\circ\f=\Phi$. Let
$$\f=U + \te1\Psi_1 + \te2\Psi_2 +\te1\te2 f$$
be the writing of $\f$ in $\mnr{n+1}$ (the algebra of
$(n+1)\times(n+1)-$matrices) and write that ${}^t\f\f=\mathbf{1}$ (it
means that if $h:=A\in \mnr{n+1}\mapsto{}^t\! AA -\mathbf{1} \in
\mnr{n+1}$, then  $\f^*(h)=h\circ\f=0$), we get
$$
\begin{array}{l}
   {}^t U U=Id\\
  A_i=U^{-1}\Psi_{i}\text{ is antisymmetric: }{}^t\! A_i=-A_i\\
  {}^t Uf +{}^t\! fU -{}^t\Psi_1\Psi_2 + {}^t\Psi_2\Psi_1=0
  \end{array}
$$
The third equation can be rewritten, setting $B=U^{-1}f$ and using
${}^t\! A_i=-A_i$,
$$
B + {}^t\! B +A_1A_2-A_2A_1=0.
$$
Now we consider the Maurer-Cartan form of $\f$:
$$
\alpha={\f}^{-1}d\f ={}^t\f d\f.
$$
We can write
$$
0=d({}^t\f\f)=(d\,{}^t\f)\f+{}^t\f d\f ={}^t\alpha +\alpha\, ,
$$
so $\alpha$ is a 1-form on $\s$ with values in $\text{so}(n+1)$.\\
Take the exterior derivative of $d\f=\f \alpha$, we get
$$
0=d(d\f)=d\f\wedge\alpha +\f d\alpha=\f(\alpha\wedge\alpha +
d\alpha).$$
Hence since $\f$ is invertible (${}^t\f\f=\mathbf{1}$)
$$
d\alpha +\alpha\wedge\alpha=0.
$$
We write $\text{so}(n+1)=\g_0 \oplus\g_1$ the Cartan decomposition
of $\text{so}(n+1)$. We have $\g_0=\text{so}(n)$ and $\g_1=\left\{
  \begin{pmatrix}
  \text{O} & v \\
    -{}^tv & 0
  \end{pmatrix},v\in \R^n \right\}$. We will write $\alpha =\alpha_0
  +\alpha_1$ the decomposition of $\alpha$.

We want to write the Euler-Lagrange equation (\ref{dbard=0}) in terms of
$\alpha$. Setting $X=\f^{-1}D\Phi$ then $\alpha_1(D)=
  \begin{pmatrix}
    \text{O} & X \\
    -{}^tX & 0
  \end{pmatrix}$ and so we have
$$
  \begin{array}{rcl}
\bar{D}X=\bar{D}(\f^{-1}D\Phi) & = &(\bar{D}\,{}^t\!\f)\f X + \f^{-1}
(\bar{D}D\Phi)\\
 & = & {}^t\alpha(D)X +\f^{-1}(\bar{D}D\Phi)
 \end{array}
$$
i.e.
\begin{equation}\label{f1}
\f^{-1}(\bar{D}D\Phi)=\bar{D}X +\alpha(\bar{D})X.
\end{equation}
Moreover
\begin{equation}\label{f2}
  \f^{-1}(\langle \bar{D}\Phi,D\Phi\rangle\Phi)=\langle\bar{D}\Phi,D\Phi
   \rangle e_{n+1}=\langle\bar{X},X \rangle e_{n+1}
\end{equation}
the last equality results from the fact that $\f$ is a map into
$\text{SO}(n+1)$; $(e_i)_{1\leq i\leq n+1}$ is the canonical basis
of $\R^{n+1}$. Besides we have
\begin{equation}\label{f3}
  \alpha(\bar{D})X=
  \begin{pmatrix}
    \alpha_0(\bar{D}) & \bar{X} \\
    -{}^t\bar{X} & 0
  \end{pmatrix}
  \begin{pmatrix}
    X \\
    0
  \end{pmatrix}=
  \begin{pmatrix}
    \alpha_0(\bar{D})X \\
    -\langle \bar{X},X\rangle
  \end{pmatrix}.
\end{equation}
Hence, combining (\ref{f1}), (\ref{f2}) and (\ref{f3}), we obtain
that the equation (\ref{dbard=0}) is written in terms of $\alpha$:
$$
\bar{D}X + \alpha_0(\bar{D})X=0\,,
$$
or equivalently
$$
\bar{D}\alpha_1(D) + [\alpha_0(\bar{D}),\alpha_1(D)]=0
$$
 where $[\,,\,]$ is the supercommutator. Thus, we have the following:

\begin{thm}
Let $\Phi\colon\s\to S^n$ be a superfield with lift $\f\colon\s\to
\text{SO}(n+1)$, then $\Phi$ is superharmonic \iif the Maurer-Cartan
form  $\alpha=\f^{-1}d\f=\alpha_0 +\alpha_1$ satisfies
$$
\bar{D}\alpha_1(D) + [\alpha_0(\bar{D}),\alpha_1(D)]=0.
$$
\end{thm}

\subsection{The general case}

We suppose that $M=G/H$ is a Riemannian symmetric space with symmetric
involution $\tau\colon G\to G$ so that $G^{\tau}\supset H\supset (G^{\tau})_0$.
Let $\pi\colon G\to M$
be the canonical projection and let $\g,\,\g_0$ be the Lie
algebras of $G$ and $H$ respectively. Write $\g=\g_0\oplus\g_1$ the
Cartan decomposition, with the commutator relations
$[\g_i,\g_j]\subset\g_{i+j\text{ mod\,2}}$.\\
Recall that the tangent bundle  $TM$ is canonically isomorphic to
the subbundle $[\g_1]$ of the trivial bundle $M\times\g$, with fiber
$\text{Ad}g(\g_1)$ over the point $x=g.H\in M$. Under this
identification the Levi-Civita connection of $M$ is just the flat
differentiation in $M\times\g$ followed by the projection on
$[\g_1]$ along $[\g_0]$ (which is defined  in the same way as
$\g_1$) (see \cite{8} and \cite{DPW}). Let $\Phi\colon\s\to M$ be a
superfield with lift $\f\colon\s\to G$ so that $\pi\circ\f=\Phi$.
Consider the Maurer-Cartan form of $\f$: $\alpha=\f^{-1}.d\f$. It is
the pullback by $\f$ of the Maurer-Cartan form of the group $G$. It
is a 1-form on $\s$ with values in the Lie algebra $\g$. We
decompose it in the form $\alpha=\alpha_0 + \alpha_1$, following the
Cartan decomposition. Then the canonical isomorphism  of bundle
between $TM$ and $[\g_1]$ leads to
a isomorphism between $\Phi^*(TM)$ and $\Phi^*[\g_1]$ and the image
of $D\Phi$ by this isomorphism is $\text{Ad}\f(\alpha_1(D))$. Thus
the Euler-Lagrange equation (\ref{dcovd=0}) is written
$$
[\bar{D}(\,\text{Ad}\f(\alpha_1(D))\,)]_{\Phi^*[\g_1]} =0
$$
where $[\cdot]_{\Phi^*[\g_1]}$ is the projection on $[\g_1]$ along
$[\g_0]$, pulled back by $\Phi$ to the projection on ${\Phi^*[\g_1]}$
along ${\Phi^*[\g_0]}$. Using the fact that
$$
A\colon (g,\eta)\in G\times\g\mapsto \text{Ad}g(\eta)$$
satisfies
$$
dA=\text{Ad}g(d\eta + [g^{-1}.dg,\eta])\,,
$$
where $g^{-1}.dg$ is the Maurer-Cartan form of $G$, this equation
becomes
\begin{eqnarray*}
0 & = & [\,\text{Ad}\f\left(\bar{D}\alpha_1(D) +
[\alpha(\bar{D}),\alpha_1(D)]\right)\,]_{\Phi^*[\g_1]}\\
 & = & \text{Ad}\f[\,\bar{D}\alpha_1(D) + [\alpha(\bar{D}),\alpha_1(D)]\,]_1\\
 & = & \text{Ad}\f\left(\bar{D}\alpha_1(D) + [\alpha_0(\bar{D}),\alpha_1(D)]
 \,\right).
\end{eqnarray*}
So we arrive at the same characterization as in the particular case
$M=S^n$.
\begin{thm}\label{theorem}
A superfield $\Phi\colon\s\to M$ with lift $\f\colon\s\to G$ is
superharmonic \iif the Maurer-Cartan form
$\alpha=\f^{-1}.d\f=\alpha_0 + \alpha_1$ satisfies
$$
\bar{D}\alpha_1(D) + [\alpha_0(\bar{D}),\alpha_1(D)]=0.
$$
\end{thm}

\section{The zero curvature equation}\label{2.4}

\begin{lemma}
Each 1-form $\alpha$ on $\s$ can be written in the form:
$$
\alpha =d\ta\,\alpha(D) + d\tb\,\alpha(\bar{D}) + (dz +
(d\ta)\ta)\,\alpha(\textstyle\dl{}{z}) +(d\bar{z} +(d\tb)\tb)
\,\alpha(\textstyle\dl{}{\bar{z}}).
$$
\end{lemma}
\textbf{Proof.} The dual basis of $\left\{D,\bar{D},\dl{}{z},\dl{}{\bar{z}}
\right\}$ is $\left\{d\ta, d\tb, dz + (d\ta)\ta , d\bar{z} + (d\tb)\tb
\right\}$.\hfill $\blacksquare$\\

\noindent We consider now that $\alpha$ is a 1-form on $\s$ with
values in the Lie algebra $\g$, then using the writing given by the
lemma, we have
\begin{equation}\label{courbnul}
\begin{array}{l}
 d\alpha + {\displaystyle\frac{1}{2}}[\alpha\wedge\alpha] = \hfill\null\\
\begin{array}{rcl}
& &  -d\ta\wedge d\ta \left\{D\alpha(D) + {\displaystyle\frac{1}{2}}[\alpha(D),
 \alpha(D)]
 + \alpha(\dl{}{z})\right\}\\
& &  -d\tb\wedge d\tb \left\{\bar{D}\alpha(\bar{D}) + {\displaystyle\frac{1}{2}}
[\alpha(\bar{D}),
 \alpha(\bar{D})] + \alpha(\dl{}{\bar{z}})\right\}\\
& & -d\ta\wedge d\tb \left\{\displaystyle\frac{}{}\!\bar{D}\alpha(D) + D\alpha(\bar{D})
  + [\alpha(\bar{D}), \alpha(D)]\frac{}{}\!\right\}\\
& &  + (dz +(d\ta)\ta)\wedge(d\bar{z}
 +(d\tb)\tb)\left\{{\displaystyle\frac{}{}}\!\partial_z\alpha(\dl{}{\bar{z}})-
 \partial_{\bar{z}}\alpha(\dl{}{z}) +
  \left[\alpha(\dl{}{z}),
 \alpha(\dl{}{\bar{z}})\right]{\displaystyle\frac{}{}}\!\right\}\\
& &  + (d\ta)\wedge(dz +(d\ta)\ta)\left\{{\displaystyle\frac{}{}}\!
  D\alpha(\!\dl{}{z})- \partial_{z}\alpha(D) + \left[\alpha(D),
 \alpha(\dl{}{z})\right]{\displaystyle\frac{}{}}\!\right\}\\
& &  + \text{ conjugate expression}{\displaystyle\frac{}{}}\\
& &  + d\ta\wedge(d\bar{z}
 +(d\tb)\tb)\left\{{\displaystyle\frac{}{}}\!D\alpha(\dl{}{\bar{z}})-
 \partial_{\bar{z}}\alpha(D) + \left[\alpha(D),
 \alpha(\dl{}{\bar{z}})\right]{\displaystyle\frac{}{}}\!\right\}\\
& &  + \text{ conjugate expression.}{\displaystyle\frac{}{}}
\end{array}
\end{array}
\end{equation}

\noindent In the following, we will write the terms like $\frac{1}{2}[\alpha(D),
 \alpha(D)]$ in the form $\alpha(D)^2$. It is justified by the fact
 that if we embedd $\g$ in a matrices algebra or more intrinsically
 in its universal enveloping algebra, so that we can write
 $[a,b]=ab-ba$, then the supercommutator is given by
 $$
 [a,b]=ab-(-1)^{p(a)p(b)}ba \,,
 $$
$p$ being the parity, and thus $[a,a]=2a^2$ if $a$ is odd.\\
The following theorem characterizes the 1-forms on $\s$ which are
Maurer-Cartan forms.

\begin{thm}\label{cbnul}\quad\\
 $\bullet$ Let $\alpha$ be a 1-form on $\s$ with values in
  the Lie algebra $\g$ of the Lie group $G$. Then there exists
  $\f\colon\s\to G$ such that $d\f=\f\alpha$ if and only if
  $$ d\alpha + \frac{1}{2}[\alpha\wedge\alpha] =0$$
  Moreover, if $U(z_0)$ is given then $\f$ is unique ($z_0\in\R^2$,
  $U=i^*\f$).\\[.2cm]
 $\bullet$ Let $A_D,\,A_{\bar{D}}\colon \s\to\g\otimes\C$ be odd
 maps, then the two following statements are equivalent
\begin{eqnarray}
   \text{(i) } & \exists \f\colon\s\to G^{\C}/D\f=\f A_D,\, \bar{D}\f=\f
   A_{\bar{D}}  & \\
    \text{(ii) } & \bar{D}A_D + DA_{\bar{D}} + [A_{\bar{D}},A_D]=0. & \label{dad}
  \end{eqnarray}
Moreover $\f$ is unique if we give ourself $U(z_0)$, and $\f$ is
with values in $G$ \iif $A_{\bar{D}}=\overline{A_D}$. In
particular, the natural map
$$
  \begin{array}{crcl}
I_{(D,\bar{D})}\colon & \negthickspace\negthickspace\negthickspace
\{\alpha\ 1\text{-form}/d\alpha + \alpha\wedge\alpha=0\} &\negthickspace
 \longrightarrow  &\negthickspace \{(A_D,A_{\bar{D}})\text{ odd which
 satisfy (ii)}\}\\
 & \negthickspace\negthickspace\negthickspace\alpha &\negthickspace
 \longmapsto &\negthickspace  (\alpha(D),\alpha(\bar{D}))
   \end{array}
$$
is a bijection.
\end{thm}
\begin{rmk}\label{remark}\emph{
$\bullet$ Suppose that $A_{\bar{D}}=\overline{A_D}$. If we embedd $\g$ in
a matrices algebra then (ii) means that:
$$
\bar{D}A_D + DA_{\bar{D}} + A_{\bar{D}}A_D + A_D A_{\bar{D}}=0
$$
i.e. $$\text{Re}(\bar{D}A_D + A_{\bar{D}}A_D)=0.$$
$\bullet$ We can see according to (\ref{courbnul}) that if $d\alpha +
\frac{1}{2}[\alpha\wedge\alpha]=0$ then $\alpha(\dl{}{z})$ (resp.
$\alpha(\dl{}{\bar{z}})$) can be expressed in terms of $\alpha(D)$
(resp. $\alpha(\bar{D})$):
\begin{equation}\label{alphaz}
\alpha\left(\dl{}{z}\right)=-(\,D\alpha(D) + \alpha(D)^2\,).
\end{equation}}
\end{rmk}

\noindent\textbf{Proof of the theorem \ref{cbnul}.} The first point follows from the Frobenius
theorem (which holds in supermanifolds, see
\cite{{2},{leites},{manin}}), for the existence. For the uniqueness, if
$\f$ and $\f'$ are solution then $d(\f'\f^{-1})=0$ so $\f'\f^{-1}$
is a constant $C\in G$, and $C=U'(z_0)U^{-1}(z_0)$.\\
For the second point, the implication (i)$\Longrightarrow$(ii)
follows from (\ref{courbnul}) (see the term in $d\ta\wedge d\tb$).
Let us prove (ii)$\Longrightarrow$(i).\\
$A_D$ and $A_{\bar{D}}$ are odd maps from $\s$ into $\g\otimes\C$ so
let us write
\begin{eqnarray*}
A_D & = & A_D^0 + \ta A_D^{\ta} + \tb A_D^{\tb} + \ta\tb
A_{D}^{\ta\tb}\\
A_{\bar{D}} & = & A_{\bar{D}}^0 + \ta A_{\bar{D}}^{\ta} +
\tb A_{\bar{D}}^{\tb} + \ta\tb A_{\bar{D}}^{\ta\tb}
\end{eqnarray*}
then we have
\begin{eqnarray*}
\bar{D}A_D & = & A_D^{\tb} - \ta A_D^{\ta\tb} - \tb\dl{A_D^0}{\bar{z}} + \ta\tb
\dl{A_D^{\ta}}{\bar{z}}\\
DA_{\bar{D}} & = & A_{\bar{D}}^{\ta} + \tb A_{\bar{D}}^{\ta\tb}
 -\ta\dl{A_{\bar{D}}^0}{z} - \ta\tb\dl{A_{\bar{D}}^{\tb}}{z}.
\end{eqnarray*}
Thus the equation (\ref{dad}) splits into 4 equations:
\begin{equation}\label{adtheta}
  \begin{array}{l}
    A_D^{\tb} + A_{\bar{D}}^{\theta} + [A_{\bar{D}}^0,A_D^0]=0  \\
\displaystyle -A_D^{\ta\tb}-\dl{A_{\bar{D}}^0}{z} + [A_{\bar{D}}^{\ta},A_D^0]
+ [A_D^{\ta},A_{\bar{D}}^0] =0\\[0.2cm]
\displaystyle A_{\bar{D}}^{\ta\tb} -\dl{A_D^0}{\bar{z}} + [A_D^{\tb},A_{\bar{D}}^0]
+ [A_{\bar{D}}^{\tb},A_D^0]=0\\
\displaystyle \dl{A_D^{\ta}}{\bar{z}} - \dl{A_{\bar{D}}^{\tb}}{z} +
[A_D^0,A_{\bar{D}}^{\ta\tb}] + [A_{\bar{D}}^0,A_D^{\ta\tb}]
+ [A_D^{\ta},A_{\bar{D}}^{\tb}] + [A_{\bar{D}}^{\ta},A_D^{\tb}]=0.
  \end{array}
\end{equation}
Now, let us embedd $\g$ in a matrices algebra $\mnr{m}$, then the Lie
bracket in $\g$ is given by $[a,b]=ab-ba$. Let us define $A,\underline{A},\,
\beta,\,B,\underline{B}$  by:
\begin{equation}\label{set}
  \begin{array}{lllllll}
 A=A_{D}^0 & , & \underline{A}=A_{\bar{D}}^0 & , & A_D^{\ta}=
- \beta(\dl{}{z}) - A^2 & , &
A_{\bar{D}}^{\tb}=-\beta(\dl{}{\bar{z}}) -
\underline{A}^2,\\[.1cm]
 & & & & A_D^{\tb} = B-\underline{A}A & , & A_{\bar{D}}^{\ta}=\underline{B}-
A\underline{A}\,,
  \end{array}
\end{equation}
then the four previous equations (\ref{adtheta}) are written:
\begin{eqnarray}
 &  B + \underline{B}=0 & \\
& A_D^{\ta\tb}={\displaystyle -\dl{\underline{A}}{z}}
+ [-B - A\underline{A},A]
+ [- \beta(\dl{}{z}) - A^2,\underline{A}] & \label{four2}\\
 &  A_{\bar{D}}^{\ta\tb} = {\displaystyle\dl{A}{\bar{z}}} +
[\underline{A}, B-\underline{A}A]
+ [A,-\beta(\dl{}{\bar{z}}) -\underline{A}^2] & \\
 & {\displaystyle\dl{}{z}}\beta(\dl{}{\bar{z}})
{\displaystyle -\dl{}{\bar{z}}}\beta(\dl{}{z})
+ {\displaystyle \dl{\underline{A}^2}{z}}
-{\displaystyle\dl{A^2}{\bar{z}}} & \nonumber\\
 & + \left[A,{\displaystyle\dl{A}{\bar{z}}}
+ [\underline{A}, B-\underline{A}A]
+ [A,-\beta(\dl{}{\bar{z}})-\underline{A}^2]
\right] & \nonumber\\
& + \left[\underline{A},{\displaystyle -\dl{\underline{A}}{z}}
+ [-B - A\underline{A},A]
+ [- \beta(\dl{}{z}) -
A^2,\underline{A}]\right] & \nonumber\\
& \,\, + \left[-\beta(\dl{}{\bar{z}})-\underline{A}^2 ,
-\beta(\dl{}{\bar{z}})-\underline{A}^2\right]
+ \left[-B-A\underline{A}, B-\underline{A}A\right]\,=\,0. &
 \end{eqnarray}
The last equation becomes after simplification
$$
{\displaystyle\dl{}{z}}\beta({\textstyle\dl{}{\bar{z}}})
 {\displaystyle -\dl{}{\bar{z}}}\beta({\textstyle\dl{}{z}})
 +[\beta({\textstyle\dl{}{z}}),\beta({\textstyle\dl{}{\bar{z}}})]=0
$$
so since $\beta$ is even and with values in $\g^{\C}$ (resp. in $\g$
if $A_{\bar{D}}=\overline{A_D}$), according to (\ref{set}) ,
we deduce from this that there exists $U\colon\s\to G^{\C}$ such that
$U^{-1}dU=\beta$ and $U$ is unique if $U(z_0)$ is given, and with values in
$G$ if $A_{\bar{D}}=\overline{A_D}$. Then we set\footnote{See remark~\ref{j}.}
\begin{equation}\label{set'}
\frac{1}{2}\Psi =UA ,\ \frac{1}{2}\underline{\Psi}=U\underline{A},\
f=\frac{2}{i} UB
\end{equation}
and
\begin{equation}\label{setf}
\f=U + \frac{1}{2}(\,\ta\Psi + \tb\underline{\Psi}\,) + \frac{i}{2}\ta\tb f .
\end{equation}
The result $\f$ is a superfield from $\s$ into $\mnc{m}$ and according to
(\ref{psitangent}) (with $\R^N=\mnc{m}$, $M=\text{GL}_m(\C)$,
$f_{\alpha}=0$, $U_{\alpha}=M$) since $U$ is invertible and hence with
values in $\text{GL}_m(\C)$, $\f$ takes values in $\text{GL}_m(\C)$.
Besides it takes values in $\text{GL}_m(\R)$ if $A_{\bar{D}}=\overline{A_D}
$. We compute that
\begin{eqnarray*}
\f^{-1} & = & \left( U + \left[\frac{1}{2}(\,\ta\Psi + \tb\underline{\Psi}\,)
+ \frac{i}{2}\ta\tb f \right]\right)^{-1}\\
 & = & \sum_{k=0}^2 (-1)^k\left[U^{-1}\left(\frac{1}{2}(\,\ta\Psi +
 \tb\underline{\Psi}\,) + \frac{i}{2}\ta\tb f \right)\right]^k
 U^{-1}\\
 & = & \left[ \mathbf{1} - (\ta A + \tb\underline{A}) - \ta\tb B +
 \ta A\tb \underline{A} + \tb\underline{A}\ta A \right]U^{-1}\\
& = & \left[ \mathbf{1} - \ta A - \tb\underline{A} - \ta\tb ( B +
  A \underline{A} - \underline{A} A ) \right]U^{-1}
\end{eqnarray*}
so
\begin{eqnarray*}
\f^{-1}.D\f & = & \f^{-1}\left(\frac{1}{2}\Psi -\ta\dl{U}{z} +
\frac{i}{2}\tb f -\ta\tb\dl{\underline{\Psi}}{z}\right)\\
 & = & A + \ta(-\beta({\textstyle\dl{}{z}}) -A^2 ) + \tb (B -
 \underline{A}A ) \\
 & & \quad +\, \ta\tb \left(-\dl{\underline{A}}{z} - \beta({\textstyle
 \dl{}{z}})\underline{A} - ( B  + A \underline{A} - \underline{A} A )A
 + AB + \underline{A}\beta({\textstyle\dl{}{z}})\right)\\
  & = & A + \ta(-\beta({\textstyle\dl{}{z}}) -A^2 ) + \tb (B -
 \underline{A}A ) \\
 & & \qquad +\, \ta\tb \left({\displaystyle -\dl{\underline{A}}{z}}
+ [-B - A\underline{A},A]
+ [- \beta({\textstyle\dl{}{z}}) - A^2,\underline{A}]\right)
\end{eqnarray*}
thus according to (\ref{set}) and (\ref{four2}) we conclude that
$$
\f^{-1}.D\f = A_D.
$$
We can check in the same way that $\f^{-1}.\bar{D}\f=A_{\bar{D}}$.
Moreover if we consider $\alpha=\f^{-1}.d\f$ the Maurer-Cartan form
of $\f$ then $(\alpha(D),\alpha(\bar{D}))=(A_D,A_{\bar{D}})$ is with
values in $\g^{\C}$, and hence it holds also for
$\alpha(\dl{}{z})$,\,$\alpha(\dl{}{\bar{z}})$ according to
(\ref{alphaz}). So $\alpha$ takes values in $\g^{\C}$. But,  according to
the first point of the theorem, the equation $\f^{-1}.d\f=\alpha$
has a unique solution if $U(z_0)$ is given, and this solution is
with values in $G^{\C}$ since $\alpha$ takes values in $\g^{\C}$
and $U(z_0)$ is in $G^{\C}$. So $\f$ takes values in $G^{\C}$ .
Moreover, $\f$ takes values in $G$ if $A_{\bar{D}}=\overline{A_D}$.
Hence, the map $I_{(D,\bar{D})}$ is surjective. Besides it is injective by
the second point of the remark~\ref{remark}: according to
(\ref{alphaz}), $\alpha$ is completely determined  by $(\alpha(D),
\alpha(\bar{D}))$. We have proved the theorem.\hfill$\blacksquare$\\
\begin{rmk}\label{j}\emph{
In general, $G$ is not embedded in $\text{GL}_m(\R)$. But since
$\g$ is embedded in $\mnr{m}$, there exists a unique morphism of group,
which is  a immersion, $j\colon G\to\text{GL}_m(\R)$, the image of
which is the subgroup generated by $\exp(\g)$. In other words $G$ is
an integral subgroup of $\text{GL}_m(\R)$ (and not a closed
subgroup). In the demonstration we use the abuse of language
consisting in identifying $G$ and $j(G)$. For example in (\ref{set'})
and  (\ref{setf}) we must use $j\circ U$ instead of $U$; and in
the end of the demonstration, when we use  the first point of
theorem, we must say that there exists a unique solution with values in
$G$, $\f_1$, and by the uniqueness of the solution (in $\text{GL}_m(\R)$)
we have $j\circ\f_1=\f$.\\
However, in the case which  interests us, $G$ is
semi-simple so it can be represented as a subgroup of $\text{GL}_m(\R)$
via the adjoint representation, and so there is no ambiguity in this
case.}
\end{rmk}
\begin{rmk}\emph{
To our knowledge, this theorem (more precisely the implication
(ii)$\Longrightarrow$(i)) has never be demonstrated in the
literature. We have only found a statement without any proof, of this
one, in \cite{odea}.}
\end{rmk}
Now we are able to prove:
\begin{thm}\label{susy}
Let $\Phi\colon\s \to M=G/H$ be a superfield into a symmetric space with
lift $\f\colon\s\to G$ and Maurer-Cartan form $\alpha=\f^{-1}.d\f$, then
the following statements are equivalent:
\begin{itemize}
  \item[(i)] $\Phi$ is superharmonic.
  \item[(ii)] Setting $\alpha(D)_{\lm}=\alpha_0(D) +
  \lm^{-1}\alpha_1(D)$ and
  $\alpha(\bar{D})_{\lm}=\overline{\alpha(D)_{\lm}}=\alpha_0(\bar{D})
  + \lm\alpha_1(\bar{D})$, we have
  $$
  \bar{D}\alpha(D)_{\lm} + D\alpha(\bar{D})_{\lm} +
  [\alpha(\bar{D})_{\lm},\alpha(D)_{\lm}]=0, \qquad \forall\lm\in
  S^1.
  $$
  \item[(iii)] There exists a lift $\f_{\lm}\colon\s\to G$ such that
  $\f_{\lm}^{-1}.D\f_{\lm}=\alpha_0(D) + \lm^{-1}\alpha_1(D)$, for
  all $\lm\in S^1$.
  \end{itemize}
Then, in this case, for all $\lm\in S^1$, $\Phi_{\lm}=\pi\circ\f_{\lm}$ is
superharmonic.
\end{thm}
\vspace{0.35cm}
\textbf{Proof.} Let us split the equation (\ref{dad}) into the sum
$\g=\g_0\oplus\g_1$:
$$
 \left\{ \begin{array}{l}
 \bar{D}\alpha_0(D) + D\alpha_0(\bar{D}) + [\alpha_0({\bar{D}),\alpha_0(D)}] +
 [\alpha_1(\bar{D}),\alpha_1(D)]=0 \\
\text{Re}(\bar{D}\alpha_1(D) + [\alpha_0(\bar{D}),\alpha_1(D)]\,)=0
  \end{array}\right.
  $$
so (ii) means that
$$
\forall \lm\in S^1,\quad \text{Re}\left(\lm^{-1}(\bar{D}\alpha_1(D) +
[\alpha_0(\bar{D}),\alpha_1(D)])\,\right)=0
$$
which means that
$$
\bar{D}\alpha_1(D) + [\alpha_0(\bar{D}),\alpha_1(D)]=0
$$
hence (i)$\iff$(ii), according to theorem~\ref{theorem}. Moreover
according to the theorem~\ref{cbnul} (ii) and (iii) are equivalent.
That completes the proof.\hfill $\blacksquare$\\[0.3cm]
We know that  the extended Maurer-Cartan form, $\alpha_{\lm}$ given by
the previous theorem is defined by $\alpha_{\lm}(D)=\alpha_0(D) +
  \lm^{-1}\alpha_1(D)$ and (so)
  $\alpha_{\lm}(\bar{D})=\alpha_0(\bar{D})  + \lm\alpha_1(\bar{D})$.
However  we want to know how the other coefficients of $\alpha$ are
transformed into coefficients of $\alpha_{\lm}$. From (\ref{alphaz})
we deduce
$$
\begin{array}{ccccc}
D\alpha_0(D) & + & \alpha_0(D)^2 + \alpha_1(D)^2 & = &
-\alpha_0(\dl{}{z})\\[0.1cm]
D\alpha_1(D) & + & [\alpha_0(D),\alpha_1(D)] &=& -\alpha_1(\dl{}{z})
\end{array}
$$
hence
$$
\begin{array}{rcccc}
(\alpha_{\lm})_0(\dl{}{z}) & = & \alpha_0(\dl{}{z}) &+&
(1-\lm^{-2})\,\alpha_1(D)^2\\[0.1cm]
(\alpha_{\lm})_1(\dl{}{z}) & = & \lm^{-1}\alpha_1(\dl{}{z}). & &
\end{array}
$$
Finally we have
\begin{multline}\label{alphalambda2}
  \alpha_{\lm}=-\lm^{-2}\alpha_1(D)^2(dz +(d\ta)\ta)
  +\lm^{-1}\alpha_1' + \alpha_0 + 2\text{Re}\left(\alpha_1(D)^2(dz
  +(d\ta)\ta)\right) \\ +\lm\alpha_1''
   -\lm^2\alpha_1(\bar{D})^2(d\bar{z} + (d\tb)\tb)
\end{multline}
where
\begin{equation}\label{where}
  \begin{array}{rcc}
\alpha_1' &=& d\ta\alpha_1(D) + (dz
+(d\ta)\ta)\,\alpha_1(\dl{}{z})\\[0.1cm]
 \alpha_1'' &=& d\tb\alpha_1(\bar{D}) + (d\bar{z} +(d\tb)\tb)\,
 \alpha_1(\dl{}{\bar{z}}).
  \end{array}
\end{equation}
So, we remark that contrary to the classical case of harmonic maps
$u\colon\R^2\to G/H$, where the extended Maurer-Cartan form is given
by $\alpha_{\lm}=\lm^{-1}\alpha_1' +\alpha_0 + \lm\alpha_1''$
(see \cite{DPW}), here  in the supersymmetric case we obtain terms
on $\lm^{-2}$ and $\lm^2$, and the term on $\lm^0$ is $\alpha_0 +
2\text{Re}\left(\alpha(D)^2(dz +(d\ta)\ta)\right)$ instead of
$\alpha_0$. Moreover, since
$\alpha_1(D)^2=\frac{1}{2}[\alpha_1(D),\alpha_1(D)]$ takes values
in $\g_0^{\C}$, we conclude that $(\alpha_{\lm})_{\lm\in S^1}$ is a
1-form on $\s$ with values in
$$
\Lm\g_{\tau}=\{\xi\colon S^1\to\g \text{ smooth}/\xi(-\lm)=\tau(\xi(\lm))\}
$$
(see \cite{DPW} or \cite{PS} for more details for loop groups and
 their Lie algebras). And so the extended lift $(\f_{\lm})_{\lm\in
 S^1}\colon\s\to \Lm G$ leads to a map $(\f_{\lm})_{\lm\in
 S^1}\colon\s\to \Lm G_{\tau}$. As in \cite{DPW}, for the classical
 case, this yields the following characterization of superharmonic
 maps $\Phi\colon\s\to G/H$.
\begin{cory}\label{extendedmc}
A map $\Phi\colon\s\to G/H$ is superharmonic \iif there exists a map
$(\f_{\lm})_{\lm\in S^1}\colon\s\to \Lm G_{\tau}$ such that
$\pi\circ\f_1=\Phi$ and
\begin{multline*}
\f_{\lm}^{-1}.d\f_{\lm}=-\lm^{-2}\alpha_1(D)^2(dz +(d\ta)\ta)
  +\lm^{-1}\alpha_1' + \tilde{\alpha}_0  +\lm\alpha_1''\\
   -\lm^2\alpha_1(\bar{D})^2(d\bar{z} + (d\tb)\tb),
\end{multline*}
where $\tilde{\alpha}_0$ and $\alpha_1$ are $\g_0^{\C}$ resp.
$\g_1^{\C}$-valued 1-forms on $\s$, and $\alpha_1',\,\alpha_1''$ are given by
(\ref{where}). Such a $(\f_{\lm})$ will be
called a extended (superharmonic) lift.
\end{cory}
\begin{rmk}\emph{
Our result for the the Maurer-Cartan form (\ref{alphalambda2}) is
different from the one obtained  in \cite{{HR1},{HR3}} or in
\cite{ki}. Because in these papers, we have a decomposition
$\g=\oplus_{i=0}^3\g_i$ with $[\g_i,\g_j]\subset\g_{i+j}$, and
$\hat{\alpha}_2$, the coefficient on $\lm^2$, is independent of
$\hat{\alpha}_1$ whereas here we have
$\hat{\alpha}_2=-\hat{\alpha}_1(D)^2(dz + (d\ta)\ta)$.
As we can see it in theorem \ref{susy}, if we decide to identify all
the Maurer-Cartan forms with their images by $I_{(D,\bar{D})}$,
$(\alpha(D),\alpha(\bar{D}))$, then  the terms on $\lm^2$ and
$\lm^{-2}$ disappear and the things are analogous to  the classical
case. In other words, it is possible to have the same formulation of
the results as for the classical case if we choose to work on
$(\alpha(D),\alpha(\bar{D}))$ instead of working on  the
Maurer-Cartan form $\alpha$. But as we will see it in the Weierstrass
representation one can not get  rid completely of the terms on $\lm^2$
and $\lm^{-2}$. So these terms are not anecdotal and constitute  an
essential difference between the supersymmetric case  and the
classical one.}
\end{rmk}
\begin{rmk}\emph{
In the following, we will simply denote by $\f$ the extended lift
$(\f_{\lm})\colon\s\to\Lm\Gt$, there is no ambiguity because we
will always precise where $\f$ takes values by writing
$\f\colon\s\to\Lm\Gt$.
Besides, given  a superharmonic map $\Phi\colon\s\to G/H$, an extended lift
$\f\colon\s\to\Lm G_{\tau}$ is determined only up to a gauge
transformation $K\colon\s\to H$ because $\f H$ is also an extended
lift for $\Phi$. Then following \cite{DPW}, we denote by
$\mathcal{SH}$ the set
$$\mathcal{SH}=\{\Phi\colon\s\to G/H \text{ superharmonic,
}i^*\Phi(0)=\pi(1)\}$$
and then  we have a bijective correspondance between $ \mathcal{SH}$
and
$$
\{\f\colon \s\to \Lm G_{\tau},\text{ extended lift, }i^*\f(0)\in H
\}/C^{\infty}(\s,H).
$$
We will note $\Phi=[\f]$.}
\end{rmk}

\section{Weierstrass-type\hfill representation\hfill  of\hfill superhar-\\ monic maps}
\label{2.5}

In this section, we shall show how we can use the method of
\cite{DPW} to obtain every superharmonic map $\Phi\colon\s\to G/H$
from Weierstrass type data.\\
We recall the following (see \cite{{DPW},{PS}}):
\begin{thm}\label{decompose}
Assume that $G$ is a compact semi-simple Lie group, $\tau \colon G\to G$ a
order $k$ automorphism of $G$  with fixed point subgroup
$G^{\tau}=H$. Let $H^{\C}=H.\mathcal{B}$ be an Iwasawa decomposition for
$H^{\C}$. Then
\begin{description}
  \item[(i)] Multiplication $\Lm G_{\tau}\times\Lm
  _{\mathcal{B}}^{+}G_{\tau}^{\C} \xrightarrow{\sim} \Lm G_{\tau}^{\C}$
   is a diffeomorphism onto.
  \item[(ii)] Multiplication $\Lm _{*}^{-}G_{\tau}^{\C}\times \Lm^{+}\Gtc
  \longrightarrow \Lm\Gtc$ is a
  diffeomorphism onto the open and dense set $\mathcal{C}=\Lm _{*}^{-}
  G_{\tau}^{\C}. \Lm^{+} \Gtc$, called the big cell.
\end{description}
\end{thm}
The above loop groups are defined by
\begin{eqnarray*}
\Lm^+\Gtc & = & \{[\lm\mapsto U_{\lm}]\in\Lm\Gtc \text{ extending
holomorphically in the unit disk}\}\\
\Lm_{\mathcal{B}}^+\Gtc & = & \{[\lm\mapsto U_{\lm}]\in\Lm^+\Gtc
/ U(0)\in \mathcal{B}\}\\
\Lm_{*}^-\Gtc & = & \{[\lm\mapsto U_{\lm}]\in\Lm\Gtc \text{ extending
holomorphically in the complement} \\
 & & \text{ \,of the unit disk and } U_{\infty}=0\}.
\end{eqnarray*}
In analogous way one defines the corresponding Lie algebras
$\Lm\gtau ,\,\Lm\gtc ,\Lm_*^-\gtc$ and $\Lm_{\mathfrak{b}}^+\gtc$
where $\mathfrak{b}$ is the Lie algebra of $\mathcal{B}$. Further we
introduce
$$
\Lm_{-2,\infty}\gtc:=\{\xi\in\Lm\gtc
/\xi_{\lm}=\sum_{k=-2}^{+\infty}\lm^k\xi_k\}.
$$
\begin{defn}\label{def}
We will say that a map $f\colon\s\to M$ is holomorphic if
$\bar{D}f=0$. We will say also that a 1-form $\mu$ on $\s$ is
holomorphic if $\mu(\bar{D})=0$ and $\bar{D}\mu(D)=0$. Moreover we
will say that $\mu$ is a holomorphic potential if $\mu$ is a
holomorphic 1-form on $\s$ with values in the Banach space
$\Lm_{-2,\infty}\gtc$ and if, writing $\mu=\sum_{k\geq -2}\lm^k\mu_k$,
we have $\mu_{-2}(D)=0$. Then noticing that a holomorphic 1-form
satisfies (\ref{dad}), we can say that the vector space
$\mathcal{SP}$ of holomorphic potentials is
$$
\mathcal{SP}={I_{(D,\bar{D})}}^{-1}\{(\mu(D),0)/\mu(D)\colon\s\to
\Lm_{-1,\infty}\gtc \text{ is odd, and }\bar{D}\mu(D)=0\}.
$$
Besides for a Maurer-Cartan form $\mu$ on $\s$ (in particular for a
holomorphic 1-form) with values in $\Lm_{-2,\infty}\gtc$ the condition
$\mu_{-2}(D)=0$ is equivalent to
$\mu_{-2}(\dl{}{z})=-(\mu_{-1}(D))^2$ according to (\ref{alphaz}).
\end{defn}
As for the classical case (see \cite{DPW}), we can construct
superharmonic maps from holomorphic potential: if $\mu\in\mathcal{SP}$
then $\mu$ satisfies (\ref{dad}), so we can integrate it
$$
g_{\mu}^{-1}.dg_{\mu}=\mu,\ i^*g(0)=1
$$
to obtain  a map $g_{\mu}\colon\s\to \Lm\Gtc$. We can decompose
$g_{\mu}$ according to theorem~\ref{decompose}
$$ g_{\mu}=\f_{\mu}h_{\mu}$$
to obtain a map $\f_{\mu}\colon\s\to \Lm \Gt$ with
$i^*\f_{\mu}(0)=1$.
\begin{thm} $\f_{\mu}\colon\s\to \Lm\Gt$ is an extended  superharmonic
lift.
\end{thm}
\vspace{0.3cm}
\textbf{Proof.} We have (forgetting the index $\mu$)
$$
\f^{-1}.d\f=\text{Ad}h(\mu)-dh.h^{-1}.
$$
But $h$ takes values in $\Lm_{\mathcal{B}}^+ \Gtc$ so that
$dh.h^{-1}$ takes values values in $\Lm_{\mathfrak{b}}^+\gtc$,
hence
$$
\left[\f^{-1}.d\f\right]_{\Lm_{*}^-\gtc}=\left[\text{Ad}h(\mu)\right]_{
\Lm_{*}^-\gtc}$$
is in the form
$$
-\lm^{-2}\alpha_1'(D)^2(dz +(d\ta)\ta)  +\lm^{-1}\alpha_1'
$$
by using the definition~\ref{def} of a holomorphic potential. But
according to the reality condition contained in the definition  of
$\Lm\gtau$:
$$\left[\f^{-1}.d\f\right]_{\Lm_{*}^+\gtc}=
\overline{\left[\f^{-1}.d\f\right]_{\Lm_{*}^-\gtc}},$$
we conclude that $\f^{-1}.d\f$ is in the same form as in the
corollary~\ref{extendedmc}, so $\f$ is an extended superharmonic
lift.\hfill $\blacksquare$\\[0.15cm]
\noindent Then according to the previous theorem we have defined a map
$$ \mathcal{SW}\colon
\mathcal{SP}\to\mathcal{SH}:\mu\mapsto[\f_{\mu}]$$
\begin{thm}\label{dbarproblem}
The map $\mathcal{SW}\colon\mathcal{SP}\to\mathcal{SH}$ is surjective
and its fibers are the orbits of the based holomorphic gauge group
$$\mathcal{G}=\{h\colon\s\to\Lm^+\Gtc,\bar{D}h=0,\, i^*h(0)=1\}$$
acting on $\mathcal{SP}$ by gauge transformations:
$$h\cdot\mu=\text{Ad}h(\mu)-dh.h^{-1}.$$
\end{thm}
\vspace{0.3cm}
\textbf{Proof.} As in \cite{DPW} it is question of solving a
$\bar{D}$-problem with right hand side in the Banach Lie algebra
$\Lm^+\gtc$:
\begin{equation}\label{dbarh}
\bar{D}h=-(\alpha_0(\bar{D}) + \lm \alpha_1(\bar{D})).h
\end{equation}
with $i^*h(0)=1$. Let us embedd $G^{\C}$ in $\text{GL}_m(\C)$ ($G$
is semi-simple). Then we set
$$h=h_0 + \ta h_{\ta} + \tb h_{\tb} + \ta\tb h_{\ta\tb}
$$
and $C=-(\alpha_0(\bar{D}) + \lm \alpha_1(\bar{D}))=C_0 + \ta C_{\ta}
+ \tb C_{\tb} + \ta\tb C_{\ta\tb}$. These are respectively writing
in $\Lm^+\mnc{m}$ and in $\Lm^+\gtc$. Then (\ref{dbarh}) splits into
\begin{eqnarray*}
h_{\tb} & = & C_0 h_0\\
-h_{\ta\tb} & = & -C_0h_{\ta} + C_{\ta}h_0\\
-\dl{h_0}{\bar{z}} & = & C_{\tb}h_0 -C_0h_{\tb}\\
\dl{h_{\ta}}{\bar{z}} &=& C_0h_{\ta\tb} + C_{\ta\tb}h_0 +
C_{\ta}h_{\tb} -C_{\tb}h_{\ta}.
\end{eqnarray*}
hence we have for $h_0$
$$
\dl{h_0}{\bar{z}}=-(C_{\tb}-C_0^2)h_0.
$$
This is a $\bar{\partial}$-problem with right hand side,
$C_0^2-C_{\ta}$, in the Banach Lie algebra $\Lm^+\gtc$ which can be
solved (see \cite{DPW}). The solutions such that $h_0(0)=1$ are
determined only up to right multiplication by elements of
$$ \mathcal{G}_0=\{h_0\colon\s\to\Lm^+\Gtc
,\,h_0(0)=1,\,\partial_{\bar{z}}h_0=0\}.$$
Then $h_{\tb}$ is given by $h_{\tb} = C_0 h_0$ so it is tangent to
$\Lm^+\Gtc$ at $h_0$. $h_{\ta\tb}$ is determined by $h_0$ and
$h_{\ta}$. So it remains to solve the equation on $h_{\ta}$ which
can be rewritten, by expressing $h_{\ta\tb}$ and $h_{\tb}$ in terms
of $h_0$ and $h_{\ta}$ in a first time, and by setting
$h_{\ta}'=h_0^{-1}h_{\ta}$ in a second time, in the following way:
$$
\dl{h_{\ta}'}{\bar{z}}=\left(\beta({\textstyle\dl{}{\bar{z}}}) +
\text{Ad}h_0^{-1}(C_0^2 -C_{\tb})\right)h_{\ta}' + \text{Ad}h_0^{-1}
\left(C_{\ta\tb} + [C_{\ta},C_0]\right)$$
where $\beta=h_0^{-1}dh_0$. Thus we obtain  an equation of the form
$$
\dl{h_{\ta}'}{\bar{z}}=a h_{\ta}' + b
$$
with $a,b\colon\R^2\to\Lm^+\gtc$, which can be solved. The solutions
such that $h_{\ta}'(0)=0$ form an affine space of which underlying
vector space is
$$\left\{h_{\ta}'\colon\R^2\to\Lm^+\gtc /\,\dl{h_{\ta}'}{\bar{z}}=
a h_{\ta}',\,h_{\ta}'(0)=0\right\}.$$
So we have solved (\ref{dbarh}). It remains to check that $h$ is
with values in $\Lm^+\Gtc$. We know that  $h_0$ takes values in
$\Lm^+\Gtc$, $h_{\ta},\,h_{\tb}$ are tangent to $\Lm^+\Gtc$ at
$h_0$. It only remains to us to check that $h_{\ta\tb}$ satisfies
equation (\ref{f''}) (or (\ref{psitangent})). But to do this we need
to know more about the embedding
$G^{\C}\hookrightarrow\text{Gl}_m(\C)$. It is possible to proceed
like that (see section~\ref{weierfield}), but we will follow another
method.\\
 Let $\gamma=dh.h^{-1}$ be the right Maurer-Cartan form of $h$.
Then by (\ref{dbarh}), we have $\gamma(\bar{D})=C$, and $C$ takes values
in $\Lm^+\gtc$, so we have to prove that $\gamma(D)$ also takes
 values in $\Lm^+\gtc$, in order to conclude that $\gamma$
 takes values in $\Lm^+\gtc$ and finally that $h$ takes values
 in $\Lm^+\Gtc$, according to the first point of the theorem~\ref{cbnul}.
 Now return to the demonstration of the theorem~\ref{cbnul},
   where we put $\gamma(D):=A_D$,
 $\gamma(\bar{D}):=A_{\bar{D}}$. Then we can see that $A_D^0,\,
 A_D^{\ta}$ take values in $\Lm^+\gtc$:
$$ A_D^0=\frac{1}{2}h_{\ta}',\
A_D^{\ta}=-\beta({\textstyle\dl{}{z}}) -(A_D^0)^2.$$
Further
\begin{eqnarray*}
A_D^{\ta\tb} & = & -\dl{A_{\bar{D}}^0}{z} + [A_{\bar{D}}^{\ta},A_D^0]
+ [A_D^{\ta},A_{\bar{D}}^0]\\
 A_D^{\tb} & = & - A_{\bar{D}}^{\theta} - [A_{\bar{D}}^0,A_D^0]
\end{eqnarray*}
according to (\ref{adtheta}); so $A_D^{\ta\tb},\, A_D^{\tb}$ are
also with values in $\Lm^+\gtc$ (these equations hold for left
Maurer-Cartan forms but we have of course analogous equations for
right Maurer-Cartan forms). Finally we have proved that $\gamma(D)$
takes values in $\Lm^+\gtc$, so we have solved (\ref{dbarh}) in
$\Lm^+\Gtc$. This completes the proof of the surjectivity (see
\cite{DPW}). For the characterization of the fibres it is the same
proof as in \cite{DPW}.\hfill $\blacksquare$\\

\noindent Let $\Phi\colon\s\to G/H$ be superharmonic with holomorphic
potential $\mu\in\mathcal{SP}$ i.e. $\Phi=[\f_{\mu}]$ where
$g=\f_{\mu}h$ and $g^{-1}.dg=\mu,\,i^*g(0)=1$. Since $g$ is
holomorphic then by using (\ref{dbar}), we can see that
$g_0=i^*g\colon\R^2\to\Lm\Gtc$ is holomorphic:
$$\partial_{\bar{z}}g_0=0.$$
Furthermore, as in \cite{DPW}, let us consider
the canonical map $\mathbf{det}\colon\Lm\Gtc\to Det^*$ (in \cite{DPW}, it
is denoted by $\tau$, see this reference for the definition of the map
$\mathbf{det}$)
and the set $|S|=(\mathbf{det}\circ g_0)^{-1}({0})$. Then
according to \cite{DPW}, since $g_0$ is holomorphic  and $\mathbf{det}\colon
\Lm\Gtc\to Det^*$ is holomorphic, then $|S|$ is discrete. But, once
more according to \cite{DPW},
$$
|S|=\{z\in\R^2 /g_0(z)\notin \text{ big cell}\}.$$
The result of this is that if we denote by $S$ the discrete set
$|S|$ endowed with the restriction to $|S|$ of the structural sheaf
of $\s$, then the restriction of $g\colon\s\to\Lm\Gtc$ to the open
submanifold of $\s$, $\s\smallsetminus S$, takes values in the big cell
(according to (\ref{psitangent}) since the big cell is a open set of
$\Lm\Gtc$). Besides using the same arguments as in \cite{DPW} we
obtain that $S\subset\s$ depends only on the superharmonic map
$\Phi\colon\s\to G/H$.
\begin{thm}
Let $\Phi\colon\s\to G/H$ be superharmonic and $S\subset\s$ as
defined above. There exists a $\g_1^{\C}$-valued  odd holomorphic
fonction $\eta$ on $\s\smallsetminus S$ so that
$$
\Phi=[\f_{\mu}]$$
on $\s\smallsetminus S$, where
$$\mu={I_{(D,\bar{D})}}^{-1}(\lm^{-1} \eta,0)=-\lm^{-2}(dz + (d\ta)\ta)\eta^2
 + \lm^{-1} d\ta\,  \eta.$$
 \end{thm}
\textbf{Proof.} It is the same proof as in \cite{DPW}.\hfill
$\blacksquare$\\

\section{ The Weierstrass representation in terms of component
fields.}\label{weierfield}

Let us consider  a map $f\colon\s\to\C^{n}$, then by using
(\ref{dbar}), $f$ is holomorphic \iif $f=u +\ta \psi$ with $u,\psi$
holomorphic on $\R^2$.\\
Further according to the definition of a holomorphic potential, we
can identify $\mathcal{SP}$ with the set of odd holomorphic maps
$\mu(D)\colon\s\to \Lm_{-1,\infty}\gtc$. Such a map is written
$$\mu(D)=\mu_D^0 + \ta\mu_D^{\ta}$$
where $\mu_D^0,\,\mu_D^{\ta}$ are holomorphic maps from
$\R^2$ into $\Lm_{-1,\infty}\gtc$,  $\mu_D^0$ being odd and
$\mu_D^{\ta}$ being even. Now, let us embedd $G^{\C}$ in
$\text{GL}_m(\C)$ so that we can work in the vector space $\mnc{m}$.
Then the holomorphic map $g\colon\s\to \Lm\Gtc$ which integrates
$$g^{-1}Dg=\mu(D),\ i^*g(0)=1$$
is the holomorphic map $g=g_0 + \ta g_{\ta}$ such that the
holomorphic  maps $(g_0,g_{\ta})$ are solution of
\begin{eqnarray*}
g_0^{-1}\dl{g_0}{z} & = & -\left(\mu_D^{\ta} + (\mu_D^0)^2\right)\\
g_0^{-1}g_{\ta} & = & \mu_D^0.
\end{eqnarray*}
Hence $g_0$ is the holomorphic map which comes from the  (even) holomorphic
potential $-(\mu_D^{\ta} + (\mu_D^0)^2)dz$ defined on $\R^2$ and
with values in $\Lm_{-2,\infty}\gtc$. So we can see that the
terms on $\lm^{-2}$ of the potential which we got rid by working on
$\mu(D)$ instead of $\mu$, reappear now when we explicit the
Weierstrass representation in terms of the component fields.\\
Remark also that $(g_0,g_{\ta})$ are the component fields of $g$. Thus we see
that the writing of a holomorphic map is the same for every embedding,
and that the third component field is equal to zero. Hence we can write
$g=g_0 + \ta g_{\ta}$ without embedding $G^{\C}$, it is at the same time
the writing of $g$ in $\Lm\Gtc$ , in $\Lm\mnc{m}$ and
for every other embedding in a vector space $\Lm\C^N$ (with $G^{\C}
\hookrightarrow \C^N$).\\
Consider, now, the decomposition $g=\f h$, and write
\begin{eqnarray*}
\f & = & U + \te1\Psi_1 +\te2\Psi_2 +\te1\te2 f\\
h & = & h_0 + \te1 h_1 +\te2 h_2 + \te1\te2 h_{12}
\end{eqnarray*}
(these are writings in $\Lm\mnc{m}$). Besides we have $g=g_0 + (\te1 +i\te2)
 g_{\ta}$. Hence we obtain
\begin{equation}\label{g=fh}
\left\{
  \begin{array}{lcl}
g_0 & = &Uh_0\\
g_{\ta} & = & \Psi_1 h_0 + Uh_1\\
i g_{\ta} & = & \Psi_2 h_0 + Uh_2\\
0 & = & Uh_{12} + fh_0 + \Psi_2h_1 - \Psi_1h_2 .
  \end{array}\right.
\end{equation}
Thus $U$ is obtained by decomposing $g_0$  which comes from a holomorphic
potential, $-(\mu_D^{\ta} + (\mu_D^0)^2)dz$, defined on $\R^2$ and
with values in $\Lm_{-2,\infty}\gtc$. So $u=i^*\Phi$ is the image
by the Weierstrass representation of this potential.\\
Then, multiplying the second and third equation of (\ref{g=fh}) by $U^{-1}$
 by the left and by $h_0^{-1}$ by the right, and remembering that
$\Lm\gtc=\Lm\gtau\oplus\Lm_{\mathfrak{b}}^+\gtc$, we obtain that
\begin{eqnarray*}
\text{Ad}h_0(\mu_D^0) & = & U^{-1}\Psi_1 + h_1h_0^{-1}\\
i\text{Ad}h_0(\mu_D^0) & = & U^{-1}\Psi_2 + h_2h_0^{-1}
\end{eqnarray*}
are the decompositions of $\text{Ad}h_0(\mu_D^0)$ resp. $i\text{Ad}h_0
(\mu_D^0)$ following the previous direct sum. In particular, we have
\begin{eqnarray}
U^{-1}\Psi_1  & = & \left[\text{Ad}h_0(\mu_D^0)\right]_{\Lm\gtau}\label{upsi1}\\
U^{-1}\Psi_2  & = & \left[i\text{Ad}h_0(\mu_D^0)\right]_{\Lm\gtau}.
\label{upsi2}
\end{eqnarray}
Finally, the third component fields $f',h_{12}'$ of $\f$ resp. $h$
are the orthogonal projections of $f$ resp. $h_{12}$ on
$U.(\Lm\gtau)$ resp. $(\Lm_{\mathfrak{b}}^+\gtc)h_0$. So by multiplying the
last equation of (\ref{g=fh}) as above and by projecting on $\Lm\gtc$ we
obtain
\begin{equation}\label{f'}
 \left[(U^{-1}\Psi_1)(h_2h_0^{-1}) - (U^{-1}\Psi_2)(h_1h_0^{-1})
 \right]_{\Lm\gtc}= U^{-1}f' + h_{12}'h_0^{-1}  .
\end{equation}
This is once again the decomposition of the left hand side
following the direct sum $\Lm\gtc=\Lm\gtau\oplus\Lm_{\mathfrak{b}}^+\gtc$.
Let us precise the orthogonal projection
$$
[\cdot]_{\Lm\gtc}\colon\Lm\mnc{m}\to\Lm\gtc .
$$
To do this it is enough to precise $[\cdot]_{\g^{\C}}\colon
\mnc{m}\to\g^{\C}$. Since $\g$ is semi-simple we can consider the
embedding
$$\text{ad}\colon \g\to\text{so}(\g)\subset\text{gl}(\g).
$$
Besides in $\text{gl}(\g)$, we have the orthogonal direct sum
$\text{gl}(\g)=\text{so}(\g)\oplus\text{Sym}(\g)$. Then for $a,b\in\text{so}
(\g)$ the decomposition of $ab$ is
$$ab=\frac{1}{2}[a,b] + \frac{ab + ba}{2}.$$
In particular for $a,b\in\g$ this decomposition is the
decomposition of $ab$ following the direct sum
$\text{gl}(\g)=\g\oplus{\g}^{\perp}$. So
\begin{equation}\label{ab}
[ab]_{\g}=\frac{1}{2}[a,b].
\end{equation}
 Now let us extend $\tau$ to
$\text{gl}(\g)$ by taking $\text{Ad}\tau$ (it is a extension because
$\tau\circ\text{adX}\circ\tau^{-1}$\hspace{0cm}$=\text{ad}(\tau(X))$).
 Then by the uniqueness of
the writing $\f  =  U + \te1\Psi_1 +\te2\Psi_2 +\te1\te2 f$ in
$\Lm\text{gl}(\g)$ and since $\Lm\text{gl}(\g)_{\tau}$ is a vector
subspace of $\Lm\text{gl}(\g)$, which contains $\Lm\Gt$, we conclude
that the previous writing is also the writing of $\f$ in
$\Lm\text{gl}(\g)_{\tau}$. So $U^{-1}f$ takes values in
$\Lm\text{gl}(\g)_{\tau}$ (and in the same way $h_{12}h_0^{-1}$ is
with values in $\Lm\text{gl}(\g^{\C})_{\tau}$\,). So, as $\tau$ commutes with
the  projection $[\cdot]_{\g^{\C}}$ (because $\tau$ preserves the
scalar product), in (\ref{f'}) it is enough to project in $\Lm\g^{\C}$
(following the direct sum $\Lm\text{gl}(\g^{\C})=\Lm\g^{\C} +
\Lm(\g^{\perp})^{\C}$)
then we automatically project in $\Lm\gtc$ (following the direct sum
$\Lm\text{gl}(\g^{\C})_{\tau}=\Lm\gtc + \Lm(\g^{\perp})_{\tau}^{\C}$\,).\\
Thus returning to the left hand side of (\ref{f'}),
this one is written
\begin{multline*}
\frac{1}{2}\left[(U^{-1}\Psi_1),(h_2h_0^{-1})\right] -
  \frac{1}{2}\left[(U^{-1}\Psi_2),(h_1h_0^{-1})\right]=\\
\frac{1}{2}\left[\left[\text{Ad}h_0(\mu_D^0)\right]_{\Lm\gtau},
\left[i\text{Ad}h_0(\mu_D^0)\right]_{\Lm^+\gtc} \right] \\-
  \frac{1}{2}\left[\left[i\text{Ad}h_0(\mu_D^0)\right]_{\Lm\gtau},
  \left[\text{Ad}h_0(\mu_D^0)\right]_{\Lm^+\gtc} \right]
\end{multline*}
by using (\ref{ab}) and (\ref{upsi1})-(\ref{upsi2}). Finally $U^{-1}f'$ is
obtained by projecting this expression on $\Lm\gtau$ following the direct sum
$\Lm\gtc=\Lm\gtau\oplus\Lm_{\mathfrak{b}}^+\gtc$. If we want $U^{-1}f$
(which depends on the embedding) we can write
$$
  (U^{-1}\Psi_1)(h_2h_0^{-1}) - (U^{-1}\Psi_2)(h_1h_0^{-1})=
  U^{-1}f + h_{12}h_0^{-1}
$$
and this is the decomposition of the left hand side following the
direct sum $\Lm\text{gl}(\g^{\C})=\Lm\text{gl}(\g)\oplus\Lm^+\text{gl}
(\g^{\C})$ (and this is also the decomposition following
$\Lm\text{gl}(\g^{\C})_{\tau}=\Lm\text{gl}(\g)_{\tau}\oplus\Lm^+\text{gl}
(\g^{\C})_{\tau}$ because all terms  of the equation are
twisted).\\[0.15cm]
Lastly, the component fields of $\Phi=\pi\circ \f_1$ are
given by: $u=\pi(U)$,$\psi_i=d\pi(U).\Psi_i$ and $F'=0$. For
example, in the case $M=S^n$, $\pi$ is just the restriction to
$\text{SO}(n+1)$ of the linear map which to a matrix associates its
last column.

\section{Primitive and Superprimitive maps with values in a
4-symmetric space.}\label{2.7}

\subsection{The classical case.}\label{debut}

Let $G$ be a compact semi-simple Lie group with Lie algebra $\g$,
$\sigma\colon G\to G$ an order four automorphism with the fixed point
subgroup $G^{\sigma}=G_0$, and the corresponding Lie algebra
$\g_0=\g^{\sigma}$. Then $G/G_0$ is a 4-symmetric space. The automorphism
$\sigma$ gives us an eigenspace decomposition of $\g^{\C}$:
$$\g^{\C}=\bigoplus_{k\in\mathbb{Z}_4}\gt_k$$
where $\gt_k$ is the $e^{ik\pi/2}$-eigenspace of $\sigma$. We have
clearly $\gt_0=\g_0^{\C}$, $\overline{\gt_k}=\gt_{-k}$ and
$[\gt_k,\gt_l]\subset\gt_{k+l}$. We define $\g_2$, $\underline{\g}_1$ and $\mathfrak{m}$
 by
$$
\gt_2=\g_2^{\C}, \
\underline{\g}_1^{\C}=\gt_{-1}\oplus\gt_1 \,\text{ and }\,
\mathfrak{m}^{\C}=\bigoplus_{k\in\mathbb{Z}_4\smallsetminus\{0\}}\gt_k ,$$
it is possible because $ \overline{\gt_2}=\gt_2$ and $\overline{\gt_{-1}}
=\gt_1$. Let us set $\g_{-1}=\gt_{-1}$,\,$\g_1=\gt_1$, $\underline{\g}_0=
\g_0\oplus \g_2$ . Then
$$\g=\underline{\g}_0\oplus\underline{\g}_1$$
is the eigenspace decomposition of the involutive automorphism
$\tau=\sigma^2$. This is also a Cartan decomposition of $\g$. Let
$H=G^{\tau}$ then $LieH=\underline{\g}_0$ and $G/H$ is a symmetric space. We
use the Killing form of $\g$ to endow $N=G/G_0$ and $M=G/H$ with a
$G$-invariant metric. For the homogeneous space $N=G/G_0$ we have
the following reductive decomposition
\begin{equation}\label{reductive}
\g=\g_0\oplus\mk
\end{equation}
($\mk$ can be written $\mk=\underline{\g}_1\oplus\g_2$) with
$[\g_0,\mk]\subset\mk$. As for the symmetric space $G/H$, we can
identify the tangent bundle $TN$ with the subbundle $[\mk]$ of the
trivial bundle $N\times\g$, with fiber $\text{Ad}g(\mk)$ over the
point $x=g.G_0\in N$. For every $\text{Ad}G_0$-invariant subspace
$\mathfrak{l}\subset\g^{\C}$, we define $[\mathfrak{l}]$ in the same
way as $[\mk]$. Then  we introduce:
\begin{defn}
$\phi\colon\R^2\to G/G_0$ is primitive if $\dl{\phi}{z}$ takes
values in $[\g_{-1}]$.
Equivalently, it means that for any lift $U$ of $\phi$, with values
in $G$, $U^{-1}\dl{U}{z}$ takes values in $\g_0\oplus\g_{-1}$.
\end{defn}
We denote by $\pi_{H}\colon G\to\ G/H$, $\pi_{G_0}\colon G\to
G/G_0$ and $p\colon G/G_0\to G/H$ the canonical projections. Let
$\phi\colon\R^2\to G/G_0$, and $U$ a lift, $\phi=\pi_{G_0}\circ U$,
and $\alpha=U^{-1}.dU$. For $\alpha$, we will use the following
decompositions:
\begin{eqnarray}
\alpha & = & \alpha_0 + \alpha_{\mk}\label{alpha1}\\
\alpha & = & \underline{\alpha}_0 + \underline{\alpha}_1
\label{alpha2}\\
\alpha & = & \alpha_2 +  \alpha_{-1} +  \alpha_0 +  \alpha_1
\label{alpha3}\\
\alpha_{\mk} & = &  \alpha_{\mk}' + \alpha_{\mk}''\label{alpha4}
\end{eqnarray}
where $\alpha_{\mk}'$ is a $(1,0)$-form and $\alpha_{\mk}''$ its
complex conjugate. Using the decomposition (\ref{reductive}), we want to write the
equation of harmonic maps $\phi\colon\R^2\to G/G_0$ in terms of the
Maurer-Cartan form $\alpha$, in the same way as for harmonic maps
$u\colon\R^2\to G/H$. Then we obtain, by using the identification
$TN\simeq[\mk]$ (and so writing the harmonic maps equation in the
form
$\left[\bar{\partial}(\text{Ad}U\alpha_{\mk}')\right]_{[\mk]}=0$):
\begin{equation}\label{harmeq}
\bar{\partial}\alpha_{\mk}' + [\alpha_0''\wedge\alpha_{\mk}'] +
[\alpha_{\mk}''\wedge\alpha_{\mk}']_{\mk}=0.
\end{equation}
Then if $[\alpha_{\mk}''\wedge\alpha_{\mk}']_{\mk}=0$, we have the same
equation as for harmonic maps into a symmetric space, and in the
same way, we can check (see \cite{12}) that the extended
Maurer-Cartan form
\begin{equation}\label{alphalm}
\alpha_{\lm} = \lm^{-1}\alpha_{\mk}' + \alpha_0 + \lm\alpha_{\mk}''
\end{equation}
satisfies the zero curvature equation
$$d\alpha_{\lm} + \frac{1}{2}[\alpha_{\lm}\wedge\alpha_{\lm}]=0.$$
Conversely, if
the extended Maurer-Cartan form satisfies the zero curvature equation and
$[\alpha_{\mk}''\wedge\alpha_{\mk}']_{\mk}=0$, then $\phi$ is
harmonic (see \cite{12}).\\
In particular if we suppose that $\phi$ is primitive then
$\alpha_{\mk}'$ takes values in $\g_{-1}$ whereas $\alpha_{\mk}''$
takes values in $\overline{\g_{-1}}=\g_1$ so  $[\alpha_{\mk}''\wedge
\alpha_{\mk}']_{\mk}=0$. Moreover let us project the Maurer-Cartan
equation for $\alpha$ onto $\g_{-1}$:
$$d\alpha_{\mk}' + [\alpha_0''\wedge\alpha_{\mk}']=0
$$ this is the harmonic maps equation (\ref{harmeq}) since
$[\alpha_{\mk}''\wedge\alpha_{\mk}']_{\mk}=0$. So a  primitive map $\phi\colon \R^2\to
 G/G_0$ is harmonic. Moreover since the extended Maurer-Cartan form
 satisfies the zero curvature equation, so we can find a harmonic
 extended lift $U_{\lm}\colon\R^2\to\Lm G$ such that
 $U_{\lm}^{-1}.dU_{\lm}=\alpha_{\lm}$. Then
 $\phi_{\lm}=\pi_{G_0}\circ U_{\lm}$ is harmonic. Besides since
 $\phi$ is primitive the decomposition
\begin{equation}\label{alpham}
\alpha=\alpha_{\mk}' + \alpha_0 + \alpha_{\mk}''
\end{equation}
is also the decomposition (\ref{alpha3}) because
$\alpha_{\mk}'\in\g_{-1}$ so
$\alpha_{\mk}'=\alpha_{-1}$,\,$\alpha_{\mk}''=\alpha_1$,\,$\alpha_2=0$
then $\alpha_{\lm}$ is a $\Lm\g_{\sigma}$-valued 1-form.
Furthermore, decomposition (\ref{alpha1}) and (\ref{alpha2}) are the
same and so the decomposition (\ref{alpham}) can be rewritten
$$
\alpha= \underline{\alpha}_1' + \underline{\alpha}_0 +
\underline{\alpha}_1''$$
 and then we can consider that $\alpha$ is the Maurer-Cartan form
associated to
$u=\pi_{H}\circ U=p\circ\phi$ with the corresponding extended
Maurer-Cartan form $\alpha_{\lm}$ given by (\ref{alphalm}). Then we
conclude that $u_{\lm}=p\circ\phi_{\lm}\colon\R^2\to G/H$ is
harmonic and $U_{\lm}$ is an extended lift for it. Moreover,
$\alpha_{\lm}$ is also a $\Lm\gtau$-valued 1-form and
$(U_{\lm})\colon\R^2\to \Lm \Gt$. So we can write that
$u=\mathcal{W}(\mu)=[U]$, where
$\mathcal{W}\colon\mathcal{P}\to\mathcal{H}$ is the Weierstrass
representation:
$$
\mathcal{W}\colon\mu\in \mathcal{P}\mapsto g \text{ \small holomorphic}
\mapsto (U,h)\in C^{\infty}(\R^2,\Lm\Gt\times\Lm_{\mathcal{B}}^+\Gtc)
\mapsto\pi_{H}\circ U_1\in \mathcal{H}
$$
between the holomorphic potentials (holomorphic 1-forms $\mu$ taking values
 in $\Lm_{-1,\infty}\gtc$) and the harmonic maps (such that
 $u(0)=H$) (see \cite{DPW}). However to obtain $\mu$ we must solve
 the following $\bar{\partial}$-problem (see \cite{DPW}):
 $$
\bar{\partial}h.h^{-1}=-(\alpha_0'' + \lm \alpha_1),$$
and since $\alpha_{\lm}$ takes values in $\Lm\gs$, this is a
$\bar{\partial}$-problem with right hand side in $\Lm^+\gsc$, so we
can find a solution $h\colon\R^2\to\Lm^+\Gsc$, $h(0)=1$.
Then the holomorphic map $g=Uh$ (it is holomorphic because $h$ is
solution of the $\bar{\partial}$-problem) takes values in
$\Lm\Gsc$ and so the potential $\mu=g^{-1}.dg$ takes values in
$\Lm\gsc$. Let us write $ \mathcal{P}_{\sigma}$ the vector subspace
of $\mathcal{P}$, of holomorphic potentials taking values in
$\Lm_{-1,\infty}\gsc=\Lm_{-1,\infty}\gtc\cap\Lm\gsc$. Then we have
proved that for each primitive map $\phi\colon\R^2\to G/G_0$ there
exists $\mu\in\mathcal{P}_{\sigma}$ such that $\phi=\pi_{G_0}\circ
U$ where $g=Uh$  and $g^{-1}.dg=\mu$. However, the decomposition
$g=Uh$ is in the same way the decomposition
$$
\Lm\Gtc\overset{\text{dec}_{\tau}}{=}\Lm\Gt .\Lm_{ \mathcal{B}}^+\Gtc    $$
but also
$$
\Lm\Gsc\overset{\text{dec}_{\sigma}}{=}\Lm\Gs .\Lm_{ \mathcal{B}_0}^+\Gsc    $$
because $g$ takes values in $\Lm\Gsc$ and because of the uniqueness of the
decomposition. We can say that the decomposition $\text{dec}_{\sigma}$ (considered
as a diffeomorphism) is the restriction of $\text{dec}_{\tau}$ to
$\Lm\Gsc$.\\
Conversely, let us prove that for any $\mu\in\mathcal{P}_{\sigma}$,
$\phi=\pi_{G_0}\circ U_{\mu}$ is primitive, so that we can conclude that
the map
$$
\mathcal{W}_{\sigma}\colon\mu\in \mathcal{P}_{\sigma}\mapsto g
\mapsto (U,h)\mapsto \phi=\pi_{G_0}\circ U_1
$$
is a surjection between $ \mathcal{P}_{\sigma}$ and the primitive
maps, i.e. that it is a Weierstrass representation for
primitive maps. So suppose that $\mu\in\mathcal{P}_{\sigma}$. Then
we integrate it: $\mu=g^{-1}.dg$, $g(0)=1$ and we decompose $g=Uh$
following $\text{dec}_{\sigma}$. Since it is also the decomposition
following $\text{dec}_{\tau}$, then we know (Weierstrass representation
$\mathcal{W}$ for the symmetric space $G/H$) that
$\alpha_{\lm}=U_{\lm}^{-1}.dU_{\lm}$ is in the form
$$
\alpha_{\lm}=\lm^{-1}\underline{\alpha}_1' + \underline{\alpha}_0 +
\lm\underline{\alpha}_1''$$
 but since $\alpha_{\lm}$ is  with values in $\Lm\gs$ (because $U$
 takes values in $\Lm\Gs$) then $\underline{\alpha}_1'\in\g_{-1}$,
 $\underline{\alpha}_0\in\g_0$, $ \underline{\alpha}_1''\in\g_1$ so
 $\phi_{\lm}=\pi_{G_0}\circ U_{\lm}$ is primitive.\\
Hence we have proved the following:

\begin{thm} We have a Weierstrass representation for primitive
maps, more precisely the map:
$$
\negthickspace \begin{array}{lccccccl}
\negthickspace\mathcal{W}_{\sigma}\colon \negthickspace &
\negthickspace\negthickspace \mathcal{P}_{\sigma}\negthickspace &
\negthickspace\xrightarrow{\mathrm{int}}
\negthickspace & \negthickspace \mathrm{H}(\C,\Lm\Gsc) &\negthickspace
 \xrightarrow{\mathrm{dec}_{\sigma}} &\negthickspace C^{\infty}(\R^2,
 \Lm\Gs\times\Lm_{\mathcal{B}_0}^+\Gsc) \negthickspace & \longrightarrow &
 \negthickspace \negthickspace \mathrm{Prim}(G/G_0)\\
\negthickspace &\negthickspace\negthickspace \mu \negthickspace
& \negthickspace\longmapsto
\negthickspace & \negthickspace g &\negthickspace \longmapsto & \negthickspace
(U,h)\negthickspace & \longmapsto &
 \negthickspace\negthickspace\phi=\pi_{G_0}\circ U_1
\end{array}
$$
is surjective. $\mathrm{H}(\C,\Lm\Gsc)$ is the set of holomorphic maps
from $\C$ to $\Lm\Gsc$, and $\mathrm{Prim}(G/G_0)$ is the set of primitive
maps $\phi\colon\R^2\to G/G_0$ so that $\phi(0)=G_0$. We can say
that $\mathcal{W}_{\sigma}$ is the restriction of the Weierstrass
representation $ \mathcal{W}$ for harmonic maps into $G/H$, to the
subspace $\mathcal{P}_{\sigma}$. More precisely, we have the
following commutatif diagram:
$$
\begin{CD}
 \mathcal{P} @> \mathrm{int}>> \mathrm{H}(\C,\Lm\Gtc) @> \mathrm{dec}_{\tau}>> C^{\infty}(\R^2,
 \Lm\Gt\times\Lm_{\mathcal{B}}^+\Gtc) @> [\pi_{H}]>> \mathcal{H}\\
 @AAA  @AAA   @AAA  @AA {[p]}A\\
  \mathcal{P}_{\sigma} @> \mathrm{int}>> \mathrm{H}(\C,\Lm\Gsc) @> \mathrm{dec}_{\sigma}>> C^{\infty}
(\R^2,\Lm\Gs\times\Lm_{\mathcal{B}_0}^+\Gsc) @> [\pi_{G_0}]>> \mathrm{Prim}(G/G_0)
\end{CD}
$$
where $[\pi_{H}](U,h)=\pi_{H}\circ U_1$, $[p](\phi)=p\circ \phi$. In
particular the image by $ \mathcal{W}$ of $ \mathcal{P}_{\sigma}$ is
the subset of $\mathcal{H}$\,: $\{u=p\circ\phi,\,\phi \text{
primitive}\}$.
\end{thm}

\subsection{The supersymmetric case.}

\begin{defn}
A superfield $\tilde{\Phi}\colon\s\to G/G_0$ is primitive if
$D\tilde{\Phi}$ takes values in $[\g_{-1}]$. Equivalently, it means
that for any lift $\f$ of $\tilde{\Phi}$, with values in $G$,
$U^{-1}.DU$ takes values in $\g_0\oplus\g_{-1}$.
\end{defn}
By proceeding as above and using the  methods we developed in the
previous sections to work in superspace, we obtain the following two
theorems:
\begin{thm}
Let $\tilde{\Phi}\colon\s\to G/G_0$ a superfield, $\f\colon\s\to G$
a lift, and $\alpha=\f^{-1}.d\f$ its Maurer-Cartan form. Then
$\tilde{\Phi}$ is superharmonic \iif
$$
\bar{D}\alpha_{\mk}(D) + [\alpha_0(\bar{D}),\alpha_{\mk}(D)] +
[\alpha_{\mk}(\bar{D}),\alpha_{\mk}(D)]_{\mk}=0.
$$
Further if $[\alpha_{\mk}(\bar{D}),\alpha_{\mk}(D)]_{\mk}=0$, then
the pair $(\alpha_0(D) + \lm^{-1}\alpha_{\mk}(D),\alpha_0(\bar{D}) +
\lm\alpha_{\mk}(\bar{D}))$ satisfies the zero curvature equation
(\ref{dad}), and so yields by $I_{(D,\bar{D})}^{-1}$ to an extended
Maurer-Cartan form $\alpha_{\lm}$. In particular, if $\tilde{\Phi}$
is superprimitive then
$[\alpha_{\mk}(\bar{D}),\alpha_{\mk}(D)]_{\mk}=0$, $\tilde{\Phi}$ is
superharmonic and $\Phi=p\circ\tilde{\Phi}\colon\s\to G/H$ is
superharmonic.
\end{thm}
\begin{thm}
We have a Weierstrass representation for superprimitive
maps, more precisely with obvious notations (according to the
foregoing):
$$
\negthickspace \begin{array}{lccccccl}
\negthickspace\mathcal{SW}_{\sigma}\colon \negthickspace &
\negthickspace\negthickspace \mathcal{SP}_{\sigma}\negthickspace &
\negthickspace\xrightarrow{\mathrm{int}}
\negthickspace & \negthickspace \mathrm{H}(\s,\Lm\Gsc) &\negthickspace
 \xrightarrow{\mathrm{dec}_{\sigma}} &\negthickspace C^{\infty}(\s,
 \Lm\Gs\times\Lm_{\mathcal{B}_0}^+\Gsc) \negthickspace & \longrightarrow &
 \negthickspace \negthickspace \mathrm{SPrim}(G/G_0)\\
\negthickspace &\negthickspace\negthickspace \mu \negthickspace
& \negthickspace\longmapsto
\negthickspace & \negthickspace g &\negthickspace \longmapsto & \negthickspace
(\f,h)\negthickspace & \longmapsto &
 \negthickspace\negthickspace\tilde{\Phi}=\pi_{G_0}\circ \f_1
\end{array}
$$
is surjective. We have the following commutatif diagram:
$$
\begin{CD}
\mathcal{SP} @> \mathrm{int}>> \mathrm{H}(\s,\Lm\Gtc) @> \mathrm{dec}_{\tau}>> C^{\infty}(\s,
 \Lm\Gt\times\Lm_{\mathcal{B}}^+\Gtc) @> [\pi_{H}]>> \mathcal{SH}\\
 @AAA  @AAA   @AAA  @AA {[p]}A\\
\mathcal{SP}_{\sigma} @> \mathrm{int}>> \mathrm{H}(\s,\Lm\Gsc) @> \mathrm{dec}_{\sigma}>>
C^{\infty}
(\s,\Lm\Gs\times\Lm_{\mathcal{B}_0}^+\Gsc) @> [\pi_{G_0}]>> \mathrm{SPrim}
(G/G_0)
\end{CD}
$$
In particular the image by $ \mathcal{SW}$ of $\mathcal{SP}_{\sigma}$
is the subset of $ \mathcal{SH}$\,:
$$\{\Phi=p\circ\tilde{\Phi},\tilde{\Phi} \text{ primitive}\}.$$
Here, the holomorphic potentials of $\mathcal{SP}_{\sigma}$
take values in $\Lm_{-2,\infty}\gsc$ and the corresponding extended
Maurer-Cartan form is in the form (\ref{alphalambda2}) but with values in
$\Lm\gs\subset\Lm\gtau$ (for example, in (\ref{alphalambda2})
$\alpha_1(D)$ takes values in $\g_{-1}$ so $\alpha_1(D)^2$ takes values in
 $[\g_{-1},\g_{-1}]\subset \g_{2}^{\C}$).
\end{thm}

\section{The second elliptic integrable system associated to a
4-symmetric space }\label{2.8}

We give us the same ingredients and notations as in the begining of
section \ref{debut}. Then let us recall what is a second elliptic system
according to C.L. Terng (see \cite{tern}).
\begin{defn}
The second $(G,\sigma)$-system is the equation for
$(u_0,u_1,u_2)\colon\C\to \oplus_{j=0}^2 \gt_{-j}$,
\begin{equation}\label{syst}
  \left\{
  \begin{array}{l}
   \partial_{\bar{z}}u_2 + [\bar{u}_0,u_2]=0\\
   \partial_{\bar{z}}u_1 + [\bar{u}_0,u_1] + [\bar{u}_1,u_2]=0\\
   -\partial_{\bar{z}}u_0 + \partial_{z}\bar{u}_0 + [u_0,\bar{u}_0] +
[u_1,\bar{u}_1] + [u_2,\bar{u}_2]=0.
  \end{array}\right.
\end{equation}
It is equivalent to say that the 1-form
\begin{equation}\label{form}
\alpha_{\lm} =\sum_{i=0}^2 \lm^{-i} u_i dz + \lm^{i} \bar{u}_i
d\bar{z}=\lm^{-2}\alpha_2' + \lm^{-1}\alpha_1' + \alpha_0 +
\lm\alpha_1'' + \lm^2\alpha_2''
\end{equation}
satisfies the zero curvature equation:
$$
d\alpha_{\lm} + \frac{1}{2}[\alpha_{\lm}\wedge\alpha_{\lm}]=0.
$$
\end{defn}
The first example of second elliptic system was given by F. H\'{e}lein
and P. Romon (see \cite{{HR1},{HR3}}): they showed that the
equations for Hamiltonian stationary surfaces in 4-dimension Hermitian
symmetric spaces are  the second elliptic system associated to
certain 4-symmetric spaces. Then we generalized the case of
$\R^4=\h$ (see \cite{HR1}) in the space $\R^8=\oct$ (with
$G=Spin(7)\ltimes\oct$, $\sigma=int_{(-L_e,0)}$, where $int_g$ is the
conjugaison by $g$, $e\in S(\text{Im}\oct)$, and $L_e$ is the left
multiplication by $e$, see \cite{ki}): there exists a family $(
\mathcal{S}_{I})$ of sets of surfaces in $\oct$, indexed by $I\varsubsetneqq
\{1,\ldots,7\}$, called the $\rho$\,-harmonic $\omega_{I}$\,-isotropic surfaces, such
that: $ \mathcal{S}_{I}\subset\mathcal{S}_{J}$
if $J\subset I$, and of which equations are the second elliptic
$(G,\sigma)$-system (see \cite{ki}).
We think that our result can be generalized to
$\oct\mathbb{P}^1,\oct\mathbb{P}^2$ or more simply to
$\h\mathbb{P}^1$.\\[0.2cm]
 For any second elliptic system associated to a 4-symmetric space,
 we can use the method of \cite{DPW} to construct a Weierstrass
 representation, defined on $ \mathcal{P}_{\sigma}^2$, the vector
 space of $\Lm_{-2,\infty}\gsc$\,-valued holomorphic 1-forms on $\C$,
 (see \cite{{HR1},{HR3}}):
$$
\mathcal{W}_{\sigma}^2\colon \mathcal{P}_{\sigma}^2\xrightarrow{\mathrm{int}}
 \mathrm{H}(\C,\Lm\Gsc) \xrightarrow{\mathrm{dec}_{\sigma}}  C^{\infty}(\R^2,
\Lm\Gs\times\Lm_{\mathcal{B}_0}^+\Gsc)   \xrightarrow{[\pi]}  \mathcal{S}
$$
where  $\mathcal{S}$ is the set of geometric maps of which equations
correspond to the second elliptic system, and $[\pi](U,h)=\pi\circ
U_1$. $\pi$ can be $\pi_{G_0}$ as well as $\pi_{H}$. For example in
the case of Hamiltonian stationary surfaces in a Hermitian symmetric
space $G/H$, we must take $\pi_{H}$ (see \cite{HR3}). Moreover if we
consider the solution $u=\mathcal{W}_{\sigma}^2(\mu)=\pi_{H}\circ
U_1$, then in this case $\phi=\pi_{G_0}\circ U_1$ can be identified
with the map $(u,e^{i\beta})$ where $\beta$ is a Lagrangian angle
function of $u$ ($G/G_0=G\times_{G_0} H$ is the principal
$U(1)$-bundle $U(G/H)/SU(2)$).
If we restrict $ \mathcal{W}_{\sigma}^2$ to $ \mathcal{P}_{\sigma}$,
we obtain $\mathcal{W}_{\sigma}$, the Weierstrass repesentation of primitive
 maps, of which image is the set of special Lagrangian surface of $G/H$
(by identifying $u$ and $\phi=(u,1)$).\\[0.25cm]
Now, we are going to give another example of second elliptic system
in the even part of a super Lie algebra. According to the previous
section, a superprimitive map $\tilde{\Phi}\colon\s\to G/G_0$ leads to
a extended lift $\f\colon\s\to\Lm\Gs$. Let
us consider $U=i^*\f\colon\R^2\to\Lm\Gs$,
then according to section~\ref{weierfield}, $U$ is obtained from
a (even) holomorphic potential, $-(\mu_{D}^{\ta} + (\mu_D^0)^2)dz$,
which is defined in $\R^2$ and with values in $\Lm_{-2,\infty}\gsc$.
 This is  a $\Lm_{-2,\infty}\gsc$-valued holomorphic 1-form on $\R^2$.
In concrete terms, if we consider that we work with the category of
supermanifolds (sets of parameters $B$, see the introduction)
$\{\R^{0|L},L\in \mathbb{N}\}$, i.e. that we work with $G^{\infty}$
functions defined on $B_{L}^{2|2}$ (see \cite{rogers}) then this is
a $(\Lm_{-2,\infty}\gsc\otimes B_L^0)$-valued holomorphic 1-form on
$\R^2$. In other words $U$ comes from a holomorphic potential which
is in $ \mathcal{P}_{\sigma}^2\otimes B_L^0$. So $u=\pi_{H}\circ
U_1\colon\R^2\to\ G/H$ as well as $\phi=\pi_{G_0}\circ
U_1\colon\R^2\to G/G_0$ correspond to  a solution of the second
elliptic system (\ref{syst}) in the Lie algebra $\g\otimes B_L^0$ (i.e.
 $u_i$ takes values in $\gt_{-i}\otimes B_L^0$). However that
 does not give us a supersymmetric interpretation of all  second
 elliptic systems (\ref{syst}) in the Lie algebra $\g$ in terms of
 superprimitive maps. Indeed, first the coefficient on $\lm^{-2}$ of
 the previous potential does not have body term: it is the square of
 a odd element so it does not have terms on $1=\eta^{\varnothing}$
(we set $B_L=\R[\eta_1,\ldots\eta_L]$). Second, this coefficient
takes values in $[\g_{-1},\g_{-1}]$ which can be $\varsubsetneqq
\g_2^{\C}$.\\[.15cm]
 In conclusion, the restrictions to $\R^2$ of
superprimitive maps $\tilde{\Phi}\colon\s\to G/G_0$ correspond to
particular solutions of the second elliptic system (\ref{syst}) in
the Lie algebra $\g\otimes B_L^0$: those which come  by
$\mathcal{W}_{\sigma}^2$, from  potentials
in the form $\hat{\mu}=-(\mu_{D}^{\ta} + (\mu_D^0)^2)dz$, with
$\mu\in\mathcal{SP}_{\sigma}$.\\
Besides for each 4-symmetric space $(G,\sigma)$, this gives us
a geometrical interpretation of certain solutions of the
second elliptic system (\ref{syst}) in $\g\otimes B_L^0$. Hence this
confirms our conjecture that there exist  geometrical problems in
$\mathbb{HP}^1$,$\mathbb{OP}^1$ and $ \mathbb{OP}^2$, analogous to
the  $\rho$\,-harmonic surfaces in $\oct$ (\cite{ki}), of which
equations are respectively the second elliptic problems in the 4-symmetric
spaces  $ \mathbb{HP}^1=Sp(2)/(Sp(1)\times Sp(1))$,
$\mathbb{OP}^1=Spin(9)/Spin(8)$ and $ \mathbb{OP}^2=F_4/Spin(9)$.

Let us give a example by considering the case of the 4-symmetric space
 $\text{SU}(3)/\text{SU(2)}$ (used by H\'{e}lein and Romon
 for their study of Hamiltonian  stationary surfaces in $\mathbb{CP}^2=
 \text{SU(3})/\text{S(U(2)}\times \text{U(1)})$).
 \begin{thm}
Consider the case of the 4-symmetric space $\mathrm{SU}(3)/\mathrm{SU}(2)$
( $H=\mathrm{S}(\mathrm{U}(2)\times \mathrm{U}(1))$). Then  an immersion
 $u\colon\R^2\to\mathbb{CP}^2(\R^{0|L})$  from $\R^2$ to the $G^{\infty}$
 manifold over
$B_L$ of $\R^{0|L}$-points of $\mathbb{CP}^2$ (morphisms from
$\R^{0|L}$  to $\mathbb{CP}^2$) is the restriction to $\R^2$ of a
superprimitive map
$$\tilde{\Phi}\colon\s\to\mathrm{SU}(3)/\mathrm{SU}(2)$$
(i.e. $u=p\circ\tilde{\Phi}\circ i$) \iif $u$ is a Lagrangian conformal
immersion
of which Lagrangian angle $\beta$ satisfies
\begin{equation}\label{beta}
\dl{\beta}{z}=ab
\end{equation}
where $a,b\colon\R^2\to\C[\eta_1,\ldots,\eta_L]$ are odd holomorphic
functions. In this case, we have $\phi=i^*\tilde{\Phi}=(u,e^{i\beta})$.
 \end{thm}
 \vspace{0.3cm}
\textbf{Proof.}
Suppose that $u$ is the restriction to $\R^2$ of a superprimitive map
$\tilde{\Phi}$, then $u$ is the image by the Weierstrass
representation $\mathcal{W}_{\sigma}^2$ of the holomorphic potential
$\hat{\mu}=-(\mu_{D}^{\ta} + (\mu_D^0)^2)dz$ with
$\mu\in\mathcal{SP}_{\sigma}$. Thus $u$ is a Lagrangian conformal immersion.
  Let us set
$$\mu_D=\lm^{-1}(A^0 + \ta A^{\ta}) + \sum_{k\geq 0}\lm^{k}\left((\mu_D^0)_k +
\ta (\mu_D^{\ta})_k\right) ,$$
 where $A^0,A^{\ta}$ takes values in $\g_{-1}$, then
$$\hat{\mu}= -\lm^{-2}(A^0)^2 dz + \sum_{k \geq -1}\lm^{k}\hat{\mu}_k .$$
Next, since $A^0$ is in $\g_{-1}\otimes B_L^1$, we can write (see
\cite{HR3})
\begin{equation}\label{a0}
A^0=\begin{pmatrix}
  0 & 0 & a \\
   0 & 0 & b \\
   -ib  & ia &  0
\end{pmatrix}
\end{equation}
thus
$$\hat{\mu}_{-2}=iab \begin{pmatrix}1&0&0 \\ 0&1&0 \\0&0& -2
\end{pmatrix}dz=3ab Ydz$$
where $Y=\frac{i}{3}Diag(1,1,-2)$. If we denote by $\hat{\alpha}_{\lm}=U^{-1}
dU=i^*\alpha_{\lm}$ the extended Maurer-Cartan form associated to $u$, then
$u$ is an immersion \iif $\hat{\alpha}_{-1}$ does not vanish.
Besides since $\g_2^{\C}=\C Y$, one can easily see  that
$$ \hat{\alpha}_2'=\hat{\mu}_{-2}$$
(because $[\g_0,\g_2]=0$). Moreover we have (see \cite{HR3})
$$\frac{d\beta}{2}Y=\hat{\alpha}_2$$
so finally
$$\dl{\beta}{z}=6ab .$$
Conversely, suppose that $u$ is a Lagrangian conformal immersion
which satisfies (\ref{beta}). Then we have $\triangle \beta=0$ since
$a,b$ are holomorphic by hypothesis. So we can write
$u=\mathcal{W}_{\sigma}^2(\hat{\mu})$ with $\hat{\mu}\in
\mathcal{P}_{\sigma}^2\otimes B_L^0$. Let us take for $\hat{\mu}$ a meromorphic
potential (see \cite{HR3})
$$\hat{\mu}=\lm^2\hat{\mu}_{-2} + \lm^{-1}\hat{\mu}_{-1}.$$
Then according to (\ref{beta}) we have $\hat{\mu}_{-2}= -(A^0)^2dz$
with $A^0$ in the same form as in (\ref{a0}). Thus if we set
$\mu_D=\lm^{-1}(A^0 - \ta\hat{\mu}_{-1}(\dl{}{z}))$, then $\mu_D$ is
an odd meromorphic map from $\s$ to $\Lm_{-1,\infty}\gsc$ and we
have
$\hat{\mu}=-(\mu_{D}^{\ta} + (\mu_D^0)^2)dz$ so
$u=p\circ\tilde{\Phi}\circ i$  with
$\tilde{\Phi}=\mathcal{SW}_{\sigma}(I_{(D,\bar{D})}^{-1}(\mu_D,0))$.\hfill
 $\blacksquare$\\

\vspace{0.2cm}
\noindent \textbf{Idrisse Khemar\\
khemar@math.jussieu.fr}


\begin{thebibliography}{111}
\bibitem{berezin} F.A. Berezin, \emph{Introduction to Superanalysis},
D.Reidel Publishing Company 1987.
\bibitem{besse} A.L. Besse, \emph{Einstein Manifolds},
Spinger-Verlag, Berlin, Heidelberg, New York, 1987.
\bibitem{12} F.E. Burstall, F. Pedit, \emph{Harmonic maps via
Adler-Kostant-Symes Theory}, Harmonic maps and integrable systems,
A.P. Fordy, J.C. Wood (Eds.), Vieweg (1994), 221-272.
\bibitem{8} F.E. Burstall and J.H. Rawnsley, \emph{Twistor theory
for Riemannian Symmetric Spaces with applications to harmonic maps
 of Riemann Surfaces} Lect. Notes in Math., vol. 1424, Springer,
 Berlin, Heidelberg, New York, 1990.
\bibitem{qfs} P. Deligne, P. Etingof, D. Freed, L. Jeffrey, D.
Kazhdan, J. Morgan, D. Morrison, E. Witten, ed. ,\emph{Quantum
Fields and Strings:  a Course for Mathematicians}, Volume 1, AMS,
1999
\bibitem{2} P. Deligne and J. Morgan, \emph{Notes on Supersymmetry},
in \cite{qfs}.
\bibitem{3} P. Deligne, D. Freed, \emph{Supersolutions}, in
\cite{qfs}.
\bibitem{DPW} J. Dorfmeister, F. Pedit and H.-Y. Wu,
\emph{Weierstrass type representation of harmonic maps into
symmetric spaces}, Comm. in Analysis and Geometry, 6(4) (1998), p.
633-668.
\bibitem{fp2} Fall Problem 2 posed by E. Witten, solutions by P.
Deligne, D. Freed, in \emph{Homework}, in \cite{qfs}.
\bibitem{Har} R. Harvey, \emph{Spinors and Calibrations},
Acadamic Press Inc., 1990.
\bibitem{HaL} R. Harvey and H. B. Lawson, \emph{Calibrated
geometries}, Acta Mathematica, 148 (1982), p. 47-157.
\bibitem{H1} F. Hélein, \emph{Applications harmoniques, lois de
conservations et repères mobiles}, Diderot éditeur, Paris 1996; or
\emph{Harmonic maps, conservation laws and moving frames}, Cambridge University
Press 2002.
\bibitem{H2} F. Hélein, \emph{Constant mean Curvature Surfaces, Harmonic
maps and Integrable Systems}, Lecture in Mathematics, ETH Zürich,
Birkhäuser 2001.
\bibitem{H3} F. Hélein, \emph{Willmore immersions and  loop groops},
J. Diff. Geometry, 50(2) (1998), p.331-338.
\bibitem{HR1} F. Hélein and P. Romon, \emph{Hamiltonian stationary Lagrangian
surfaces in $\C^2$}, Comm. in Analysis and Geometry Vol. 10, N. 1, 2002, p.
79-126.
\bibitem{HR2} F. Hélein and P. Romon, \emph{Weierstrass representation
of Lagrangian surfaces in four dimensional spaces using spinors and
quaternions}, Comment. Math. Helv., 75 (2000), p. 668-680.
\bibitem{HR3} F. Hélein and P. Romon, \emph{Hamiltonian stationnary
Lagrangian surfaces in Hermitian symmetric spaces}, in
\emph{Differential Geometry and Integrable Systems}, Martin Guest,
Reiko Miyaoka, and Yoshihiro Ohnita, Editors-AMS, 2002.
\bibitem{Hel} S. Helgason, \emph{Differential geometry, Lie group
and symmetric spaces}, Academic Press, Inc., 1978.
\bibitem{ki} I. Khemar, \emph{Surfaces isotropes de $\oct$ et
syst\`{e}mes int\'{e}grables.}, preprint arXiv:math.DG/0511258.
\bibitem{leites} D.A. Leites, \emph{Introduction to the theory of
Supermanifolds}, Russian Math. Surveys, \textbf{35} no. 1 (1980),
1-64.
\bibitem{manin} Y.I. Manin, \emph{Gauge field theory and Complex
Geometry}, Gundlehren der Mathematischen Wissenshaften 289,
Springer-Verlag, 1988.
\bibitem{odea} F. O'Dea {Supersymmetric Harmonic Maps into Lie
Groups}, preprint arXiv:hep-th/0112091.
\bibitem{PS} A. Pressley and G. Segal, \emph{Loop groops}, Oxford
Mathematical Monographs, Clarendon Press, Oxford, 1986.
\bibitem{rogers} A. Rogers, \emph{A global theory of
supermanifolds}, J.Math.Phys. \textbf{21} (6), 1980, 1352-65;
\emph{Super Lie groups: global topology and local structure},
J.Math.Phys. \textbf{22}(5), 1981, 939-45.
\bibitem{tern} C.L. Terng, \emph{Geometries and Symmetries of Soliton
Equations and Integrable Elliptic Equations}, preprint
arXiv:math.DG/0212372.
\end{thebibliography}
\end{document}